\documentclass[11pt,reqno]{amsart}

\usepackage{graphicx} 
\usepackage{bm} 

\usepackage[utf8]{inputenc}
\usepackage[T1]{fontenc}
\usepackage{mathptmx}
\usepackage{etoolbox}
\usepackage[margin=2cm]{geometry}

\usepackage{color}
\usepackage{hyperref}
\usepackage{soul}

\usepackage[dvipsnames]{xcolor}

\usepackage{subcaption}

\usepackage{tikz}
\usetikzlibrary{arrows.meta,positioning}

\usepackage{pgfplots}
\usepackage{pgfplotstable}
\usepgfplotslibrary{groupplots}
 \pgfplotsset{compat=1.18}

\usepackage{siunitx}

\usepackage{enumitem}

\usepackage{tabularray}
\UseTblrLibrary{diagbox}

\colorlet{mylinkcolor}{blue}
\colorlet{mycitecolor}{ForestGreen}
\colorlet{myurlcolor}{Aquamarine}

\hypersetup{
  linkcolor  = mylinkcolor!70!black,
  citecolor  = mycitecolor,
  urlcolor   = myurlcolor,
  colorlinks = true,
}

\usepackage{bbm}

\newcommand{\real}{\mathrm{Re}}
\newcommand{\R}{\mathbb{R}}
\newcommand{\N}{\mathbb{N}}
\newcommand{\xd}{\,\textrm{d}}
\DeclareMathOperator{\midint}{mid}
\DeclareMathOperator{\conv}{Conv}
\DeclareMathOperator{\Span}{span}
\DeclareMathOperator{\diag}{diag}

\newcommand{\sumi}{\sum_{i=1}^N}

\newcommand{\oinfmicro}{\omega_\infty^\textrm{micro}}
\newcommand{\oinfcont}{\omega_\infty^\textrm{cont}}

\newtheorem{lemma}{Lemma}[section]

\newtheorem{proposition}[lemma]{Proposition}

\newtheorem{assumption}[lemma]{Assumption}
\newtheorem{remark}[lemma]{Remark}

\usepackage{orcidlink}

\keywords{opinion formation, network mediated dynamics, systems of nonlinear conservation laws, large population limits}
\subjclass[2020]{
35Q91, 
35L65, 
35F20, 
82B40, 
91D30.
}

\begin{document}
\title[A framework for continuum modeling of opinion dynamics on a network]{A framework for continuum modeling of opinion dynamics on a network based on probability of connections}

\author{G.~Favre\,\orcidlink{0000-0002-2090-2648}}
\address{Faculty of Mathematics, University of Vienna, Oskar-Morgenstern-Platz 1, 1090, Vienna (Austria).}
\email{gianluca.favre@univie.ac.at}

\author{G.~Jankowiak\,\orcidlink{0000-0002-9025-1465}}
\address{Department of Mathematics and Scientific Computing, University of Graz, Heinrichstraße 36, 8010, Graz (Austria).}%
\email{gaspard.jankowiak@uni-graz.at}

\author{S.~Merino-Aceituno\,\orcidlink{0000-0003-2114-9210}}%
\address{Faculty of Mathematics, University of Vienna, Oskar-Morgenstern-Platz 1, 1090, Vienna (Austria).}
\address{Complexity Science Hub (External Faculty member), Josefst\"adter Strasse 39, 1080, Vienna (Austria).}
\email{sara.merino@univie.ac.at}

\author{L.~Trussardi\,\orcidlink{0000-0001-6255-9983}}
\address{Institute of Mathematics and Scientific Computing, University of Graz, Heinrichstraße 36, 8010 Graz (Austria)}
\email{lara.trussardi@uni-graz.at}
\date{\today}

\begin{abstract}

We propose a modeling framework to develop a continuum description of opinion dynamics on networks as an alternative to other models, like the ones based on graphons. In a nutshell, the continuum model that we propose aims at approximating the distribution of opinions as well as the probability that two given opinions are connected.  To illustrate our framework, we focus on a simple model of consensus dynamics on a network and derive a continuum description using techniques inspired by mean-field limits. We also discuss the limitations of this approach and suggest extensions to account for dynamic networks with evolving connections, stochastic effects, and directional interactions.
\end{abstract}

\maketitle
\thispagestyle{empty}

\section{Introduction}

The continuum description of graphs is a challenging mathematical problem, particularly the derivation of continuum models from discrete ones. A well-established approach is through graphons~\cite{Kuehn-21,Kuehn-23,Kuehn-24}, typically for dense graphs, and their extensions, such as graphops~\cite{Kuehn-20}, for sparse graphs. As an alternative to continuum descriptions of networks based on graphons, we propose here a description based on two quantities: the opinion density $f=f(\omega)$, and the edge distribution $g=g(\omega, m)$. The density $f$ represents the proportion of individuals with opinion $\omega$, and it is the limiting object obtained through classical mean-field limits. In the discrete setting, $f(\omega)$ can be understood as the number of individuals with opinion $\omega$ divided by the total population. The edge distribution $g$ gives the proportion of links between pairs of individuals with opinions $\omega$ and $m$. At the discrete level, $g(\omega,m)$ can be understood as the number of connected individuals with opinions $m$ and $\omega$ divided by the total number of connections in the network. The edge distribution $g$ can be interpreted as the probability of two individuals holding different opinions being connected, but it does not encode the full structure of the discrete network. We propose a system of partial differential equations governing the evolution of $f$ and $g$ by applying mean-field techniques. The primary challenge is that unlike classical systems, network dynamics cannot be fully captured by averaged behavior alone. By adding an extra piece of information - the edge distribution $g$, - we aim to bridge this gap.

Mean-field limits have already been successfully employed to study network dynamics. In~\cite{degond2016continuum} for instance, the authors derived continuum equations for networks of fibers, with further developments in~\cite{barre2018particle,barre2017kinetic,barre2021fast}. In that model, fibers are characterized by spatial position and orientation. In contrast, opinion dynamics involves a single variable: the opinion. This fundamental difference limits the applicability of the approach in~\cite{degond2016continuum} to opinion models, as we will illustrate later. Rigorous mean-field limits have been proven rigorously in, for example, \cite{AD24,JPS21,paul2022microscopic}. In \cite{JPS21}, the authors prove rigorously the mean-field limit under very mild assumptions on the network and obtain an equation for $h=h(t,x,\xi)$, where $x\in\mathbb{R}^d$ is the position variable and $\xi$ is an element of the network $\xi\in[0,1]$. They obtain an equation for $h$ that is computed using the extended graphon associated to the network. In this case $h=h(t,x,\xi)$ describes the probability at time $t>0$ of finding an agent at position $x\in\mathbb{R}^d$, in a certain state of interaction with the network $\xi$ - where $\xi$ would be the analog of a particle label. How the function $h$ connects with the functions $f$ and $g$ described above is unclear, as the two approaches offer different ways of describing the network - $h$ gives the connectivity of one individual, while $(f,g)$ gives the connectivity of pairs of individuals. This type of description has been also used in other mean-field limits with some variations in  \cite{APD21,paul2022microscopic}, see also the review \cite{AD24}.

\bigskip

To illustrate the difference between mean-field dynamics and network-induced dynamics, consider a classical mean-field set-up where an agent $i$ with opinion $\omega$ interacts with all agents $j$ with opinion $m$ if $|\omega-m|< R$ for some $R>0$. Crucially, determining if $i$ and $j$ interact only requires knowledge of the respective opinions $\omega$ and $m$, not their individual labels. This enables a mean-field description that is label-independent and based on the concept of exchangeability: interchanging labels yields statistically equivalent systems.

In network-based dynamics the situation is different. Whether individuals with opinions $\omega$ and $m$ influence each other depends on the connectivity between them, and such connectivity is described through the individuals' labels. Two individuals with opinions $\omega$ and $m$ may be connected, while another pair with the same opinions may not. Thus, in network-mediated systems, the labels and their associated connections are critical, making it impossible to determine if two individuals interact from just knowing their individual opinions.

To perform a mean-field limit for the opinion model, we need a label-independent representation of the network. For this, instead of tracking an individual's label and opinion, $(i, \omega)$, we consider only its opinions $\omega$. If an individual $i$ with opinion $\omega$ is connected to an individual $j$ with opinion $m$, we just say that $\omega$ is connected to $m$, disregarding labels. This approach works in simple cases but fails when multiple individuals share the same opinion $\omega$. We will discuss strategies to address this limitation, which involve augmenting the system with additional information to distinguish between individuals with identical opinions. Each specific system will require a tailored approach to incorporate this information. As an example, we demonstrate how this can be achieved for the Lancichinetti-Fortunato-Radicchi (LFR) type networks, which mimic community structures.

\bigskip

The continuum system of equations for $f$ and $g$ that we propose is formally derived under several assumptions. To test the accuracy of the continuum model, first, we verify that it preserves key structural properties of the discrete system; second, we compare the discrete and continuum models numerically.

\subsection{The opinion model and the associated network.}

Our goal is to introduce a continuum modeling framework for opinion dynamics on a network. To this end, we consider the Abelson model~\cite{abelson_math_models_1967}, a simple opinion consensus model where individuals interact pairwise through their social network. At the discrete level, we consider $N>1$ individuals, each characterized by an opinion $\omega_i(t) \in \mathbb{R}$ at time $t \in \mathbb{R}_+$ with $i=1,\dots, N$, connected via a fixed network.
The key difference between this model and the Hegselmann-Krause (HK) model~\cite{Hegselmann-02} lies in the interaction mechanism. In the HK model, agents interact based on the proximity of their opinions, whereas in Abelson's model, interactions occur only between agents connected by a pre-existing, fixed network. We assume that the network is undirected, a feature observed in social networks like BeReal, Snapchat, and earlier versions of Facebook.

For the numerical comparison between the discrete and continuum models, we use LFR-type graphs~\cite{LFR}, a well-known benchmark for testing community behavior.  The population is divided into a fixed number of communities, where individuals within the same community share similar opinions. We perform simulations by varying the number of connections between communities, ranging from well-mixed to segregated networks.

In general, consensus is eventually reached as long as the network does not have disconnected components. However, the rate of convergence to consensus depends on the network structure and the connectivity between agents. We investigate this property through numerical simulations for both the discrete and continuum models.

\subsection{Some background on social models and social networks}

The first formal models of the dynamics of opinion formation in a group were inspired by~\cite{Abelson64, French-56, Harary59}. The basic idea in these works is that if two members of a group interact, "each member of the group changes his attitude position towards the other by some constant fraction of the ‘distance’ between them"~\cite{Abelson64}. We note thus that we have a link with the concept of homophily.
In the works of~\cite{Abelson64,Berger-81,DeGroot-74,Harary59,Lehrer75}, it was analytically shown that for a broad class of models with repeated social influence, consensus is always reached unless the network of interactions is disconnected. More complex models models exist that consider different types of networks and interactions leading to phenomena like polarization \cite{SLSS23}. For a very good overview of models of social influence, we suggest the reading of~\cite{Models17}.

In these models, individuals are treated as nodes in a network. Nodes are connected by a network link if they exert influence on each other. The network links are assumed to be fixed, but they could be weighted. The weights represent the strength of social influence that one individual exerts upon another. This could represent, for example, the social status, or the frequency of interaction. Models as~\cite{Hegselmann-02} of assimilative influence, with continuous opinions and fixed influence weights, are also often called "classical" models in the literature. Such classical averaging models assume that opinions vary on a continuous scale and implement social influence as averaging~\cite{Friedkin-90, friedkin11}.

The model we suggest here falls into such a category. The key idea is that individuals which are connected always influence each other towards reducing opinion differences. As a result, if the network is fully connected, the dynamics (with continuous opinions) inevitably generate consensus in the long run.

In 2000, a sociophysic model for opinion dynamics~\cite{Sz00} including ``social validation'' and ``discord destroys'' effects was studied. Afterward, in~\cite{slanina2003analytical}, such a model was considered on a network.

Another important branch of models has been inherited from the study of birds flocking~\cite{CK2,CK1} --- or more generally from the study of the behavior of groups of animals --- introducing bounded confidence models. The work of Deffuant and Weisbuch~\cite{Deffuant-00} is one of the most known examples: it presents randomness in the choice of the interactions pairs, which are randomly chosen at each time step.

Several properties like convergence to consensus, clustering, dependence of the equilibrium from the initial data have been well investigated and studied (see for example~\cite{Berger-81,Ding-19,Dong-17,JM14,Tian-21}).
The role of the network's structure in the last decade is gaining more and more importance since, nowadays, opinion travels mainly through social media and social networks. Different models for opinion formation, also including a network structure, attracted a lot of interests~\cite{Berner-23,Salvarani21,Brooks-20,Burger-22,Favre-23,Fagioli-20,Kuehn-23,KS24,Motsch-14,Nugent-23,NGW24,Sawicki-23,Stokes-24}.
Anyhow, the description of opinion on networks covers a wide range of fields, from mathematics~\cite{Albi-16,Burger-21,Burger-22,Kuehn-21,Tosin-22}, to philosophy~\cite{Dalege-23,Rosenstock17,Zollmann-12}, through computational social science~\cite{Starnini-20}.

As already mentioned in the introduction, there is a growing interest in the kinetic theory community to investigate group behavior in social systems, and a correspondingly growing number of contributions~ \cite{degond2017continuum,dimarco2021kinetic,Bertram1,during2009boltzmann,Bertram2,jabin2023mean,toscani2018opinion} (just to name a few). The passage from discrete to continuum descriptions has been investigated in \cite{dubovskaya22,APD21,AD24,GK22,JPS21,paul2022microscopic}.

\subsection{Structure of the paper}
In Section~\ref{sec:iba}, we describe the discrete model. In Section~\ref{sec:consensus_discrete} we consider it in the specific case of consensus dynamics and offer a review of mathematical properties of convergence to consensus. In Section~\ref{sec:mfl}, we propose a continuum model to approximate the discrete dynamics when the number of individuals grows large. We start by defining the node and edge distributions in Section~\ref{sec:network description}. After the derivation of the continuum equations in Sections~\ref{sec:IBM_nolabels} and~\ref{sec:limit}, we investigate its mathematical properties in the case of consensus dynamics in Section~\ref{sec:consensus continuum}. In Section~\ref{sec:refinements}, we present refinements of the continuous model. In Section~\ref{sec:numerics}, we compare numerically discrete and continuum models.
Finally, we present extensions of the model in Section~\ref{sec:extensions_model}. Limitations and open problems are discussed in Section~\ref{sec:conclusions}.

\section{Discrete dynamics}\label{sec:iba}
\subsection{General model}
We consider a  population of $N> 1$ individuals, each having a certain opinion that evolves over time: $\omega_i=\omega_i(t)$, $\omega_i:\R^+\rightarrow\Omega$, with $i=1,\dots,N$ and $\Omega\subset \mathbb
R$ a connected set.

Moreover, we assume that individuals are part of a network, which is modeled as a unweighted, undirected graph, where each node corresponds to an individual.
Such a graph can be described by the \emph{adjacency matrix} $A=(A_{ij})\in\{0,1\}^{N\times N}$.
We say that two individuals $i$ and $j$ are connected if there is an edge between the corresponding nodes, and we then write $i \sim j$.
In terms of the entries $A_{ij}$, we have $A_{ij} = 1$ if and only if $i \sim j$.
In particular, we will assume the following.

\begin{assumption} \label{as:connected_graph}
    The graph $A$ is connected, and $A_{ij}=A_{ji}$ and $A_{ii}=0$ for all $i,j=1,\hdots,N$. This assumption means that there is no self-interaction.
\end{assumption}

The evolution of each opinion $\omega_i=\omega_i(t)$ is given by:
\begin{equation}\label{eq:ODE_eb}
    \frac{\xd}{\xd t} \omega_i(t) = \frac{1}{\# \mathcal{I}_i} \sum_{j\in \mathcal{I}_i} \mathbf{D}(\omega_i(t) - \omega_j(t)).
\end{equation}
The function $\mathbf{D} : \R \to \R$ is the so-called \emph{debate operator} and models pair-wise interactions.
The set $\mathcal{I}_i$, defined as $\mathcal{I}_i = \{j : A_{ij} = 1\}$, represents the individuals $j$ connected with the individual $i$ through the network; and $\#\mathcal{I}_i$ denotes the cardinality of this set.   The presence of the weight $(\#\mathcal{I}_i)^{-1}$ appears also in other models for collective dynamics and opinion formation~\cite{degond2017continuum,Motsch-14}  and it models that when two individuals interact, the one with more connections gets less
influenced than the other
(see e.g.~\cite{Tosin-24}).
We also note that the model can be extended to $\Omega\subset \R^d$, for $d\geq 1$, but here we will consider only the one-dimensional case.

\subsection{Case of consensus dynamics}
\label{sec:consensus_discrete}

To compare the discrete model~\eqref{eq:ODE_eb} and the continuum model proposed in this article (Section~\ref{sec:mfl}), we will focus on the particular case of consensus dynamics, for which the system reaches equilibrium asymptotically in time. We will compare both models in terms of their asymptotic behavior.

We assume that $\Omega=(-1,1)$
and that during pair-interactions  individuals' opinions relax towards a common opinion (consensus). To capture this consensus dynamics, we make  the following  assumption:

\begin{assumption}\label{ass:debate operator}
  It holds that $\mathbf{D}(z)=-\partial_z W(z)$ for all $z \in \R$, where $W\geq 0$ is a smooth strictly convex potential with a unique global minimum at $z=0$. Moreover, we choose $W$ even, i.e. $W(z)=W(|z|)$.
\end{assumption}

Note that this implies that $\mathbf{D}$ is decreasing, with $\mathbf{D}(0) = 0$. The simplest example is the quadratic potential $W(z)=z^2/2$, and in this case $\mathbf{D}(\omega_i-\omega_j) = \omega_j-\omega_i$.\\

Through the whole paper, we are going to work under Assumptions~\ref{as:connected_graph} and~\ref{ass:debate operator}. These are enough to ensure global existence and uniqueness of the solutions.

\subsubsection{Convergence to equilibrium}
In this Section, we review some mathematical properties of~\eqref{eq:ODE_eb}. In particular, we focus on the convergence to consensus. At the particle level, the following holds:
\begin{lemma}[Conserved quantity]
 \label{lem:conserved_quantity_discrete} For all times, it holds that
  \begin{align} \label{eq:dt}
    \sum_{i=1}^N \#\mathcal{I}_i \omega_i(t)=\sum_{i=1}^N \#\mathcal{I}_i \omega_i(0).
 \end{align}
\end{lemma}

\begin{proof}
\label{appendix: proof cons quant}
We compute the time derivative of the left-hand-side of \eqref{eq:dt}, which gives
\begin{align*}
    \sumi \#\mathcal{I}_i \Bigl(\frac{1}{\#\mathcal{I}_i}\sum_{j\in\mathcal{I}_i}  \mathbf{D}(\omega_i - \omega_j)\Bigr)
    = \sumi \sum_{j\in\mathcal{I}_i} \mathbf{D}(\omega_i - \omega_j)=0\,.
\end{align*}

To prove that the last expression is equal to zero, we recall now that $\mathbf{D}$ is odd, i.e. $\mathbf{D}(\omega_j - \omega_i) = - \mathbf{D}(\omega_i - \omega_j)$. Additionally, we observe that if an individual $i$ is connected with an individual $j$, i.e. $j\in\mathcal{I}_i$, then it also holds $i\in\mathcal{I}_j$. These two properties imply that all the terms in the sum cancel out, and the sum is thus equal to zero.
\end{proof}

Lemma~\ref{lem:conserved_quantity_discrete} shows that the total weighted opinion  of the network \eqref{eq:dt} is a conserved quantity.  We emphasize that Assumption~\ref{ass:debate operator} is fundamental for this conservation property to hold.

\begin{proposition}[Convergence to consensus]\label{pr:long_time_behaviour}
 Suppose that the graph is connected (Assumption~\ref{as:connected_graph}), then it follows that
$$\omega_i(t) \to \oinfmicro  \mbox{ as } t\to\infty,\forall i,$$
 where $\oinfmicro$ is the weighted averaged opinion
\begin{align} \label{eq:limit_omega}
 \oinfmicro := \frac{1}{\sum^N_{j=1} \#\mathcal I_j}\sum^N_{i=1} \#\mathcal I_i \omega_i(t=0).
\end{align}
\end{proposition}

\begin{remark} We observe that the connectedness Assumption~\ref{as:connected_graph} is essential, without it one would expect different limiting values in general, one for each of the connected components.
\end{remark}

The proof of Proposition~\ref{pr:long_time_behaviour} is based on a Lyapunov functional.

\begin{proof}[Proof of Prop. \ref{pr:long_time_behaviour}]

Define the global potential
\begin{equation} \label{eq:potential_V}
V(\omega_1,\hdots,\omega_N):= \frac{1}{2}\sum_{i=1}^N\sum_{j\in \mathcal{I}_i} W(\omega_i-\omega_j).
\end{equation}
Due to the Assumption~\ref{ass:debate operator}, it holds that
$$
\frac{\xd \omega_i(t)}{\xd t} = \frac{1}{\# \mathcal I_i}\sum_{j\in \mathcal{I}_i}\mathbf{D}(\omega_i-\omega_j) =-\frac{1}{\# \mathcal I_i}\sum_{j\in \mathcal{I}_i}\partial_{\omega_i}W(\omega_i-\omega_j) =-\frac{1}{\# \mathcal{I}_i}\partial_{\omega_i}V,
$$
where we used that $W$ is even to compute $\partial_{\omega_i}V$.

Thus
\begin{equation}
    \frac{d}{dt}V(\omega_1(t),\hdots, \omega_N(t)) = \sum_{i=1}^N\partial_{\omega_i}V\,\frac{\xd\omega_i}{\xd t} =-\sum_{i=1}^N\frac{1}{\# \mathcal I_i}(\partial_{\omega_i}V)^2\leq - \min\{(\#\mathcal I_1)^{-1},\hdots, (\#\mathcal I_N)^{-1}\}|\nabla V|^2 \leq 0. \label{eq:Lyapunov_condition_micro}
\end{equation}
 Since $V$ is the sum of strictly convex functions, it is strictly convex and in this case it has a unique critical point corresponding to the minimum. Indeed, the minimum of $V$ is attained when
$$W(\omega_i-\omega_j) =0 \qquad \mbox{for all } i\sim j,
$$
which holds if and only if $\omega_i=\omega_j$ for all $i\sim j$, by Assumption \ref{ass:debate operator}. Since we assume that the graph has only one connected component, this implies that all constant solutions, i.e. $\omega_i=\oinfmicro$ for all $i$ for some $\oinfmicro$ are minimizers of $V$. At the same time, thanks to having a conserved quantity (Lemma~\ref{lem:conserved_quantity_discrete}), we know that it must hold that
$$\sum_{i=1}^N \#\mathcal I_i \omega_i(0)= \sum_{i=1}^N \#\mathcal I_i \oinfmicro,$$
which determines the value of $\oinfmicro$ given in \eqref{eq:limit_omega}.
 This means that $V$ is a global Lyapunov function under the dynamics (thanks to \eqref{eq:Lyapunov_condition_micro}) for $\omega_i=\oinfmicro$. We can apply La Salle's principle~\cite{LaSalle} and conclude that this $(\oinfmicro)_{i=1}^N$ is a global attractor.
\end{proof}

\subsubsection{Speed of convergence in the linear case.}
  \label{sec:exponential convergence}
We will use the simplest scenario corresponding to $\mathbf D(z)=-z$ to carry out the numerical simulations in Section~\ref{sec:numerics}.
In this case, it turns out that $\omega$ converges to $\oinfmicro$ \emph{exponentially fast}. More precisely we have the following
\begin{proposition}
  \label{prop:expoconv}
If $\mathbf D(z) = -z$, there exist positive constants $C$ and $\alpha$ such that
\begin{equation}
  \label{eq:definition exponential convergence}
  \|\omega_i(t) - \oinfmicro \| \le C \exp\left(- \alpha t\right) \|\omega_i(0) - \oinfmicro\|.
\end{equation}
\end{proposition}
In the present case, this is a standard result of the theory of linear systems, a proof is sketched in Appendix~\ref{sec:app:expoconv}.
This result also holds in more general cases (e.g. time-dependent graphs), see the review of Proskurnikov and Tempo and the references therein~\cite{Proskurnikov_2018}.

\section{Continuum modeling framework: the basics}

\label{sec:mfl}

When the number of agents is very large (like in many social networks), we are interested in understanding the evolution of the opinions at the level of the population rather than tracking the changes of each individual opinion. For this purpose, a continuum description is more appropriate than a discrete one.

\begin{figure}
    \centering
    \includegraphics[width=\linewidth]{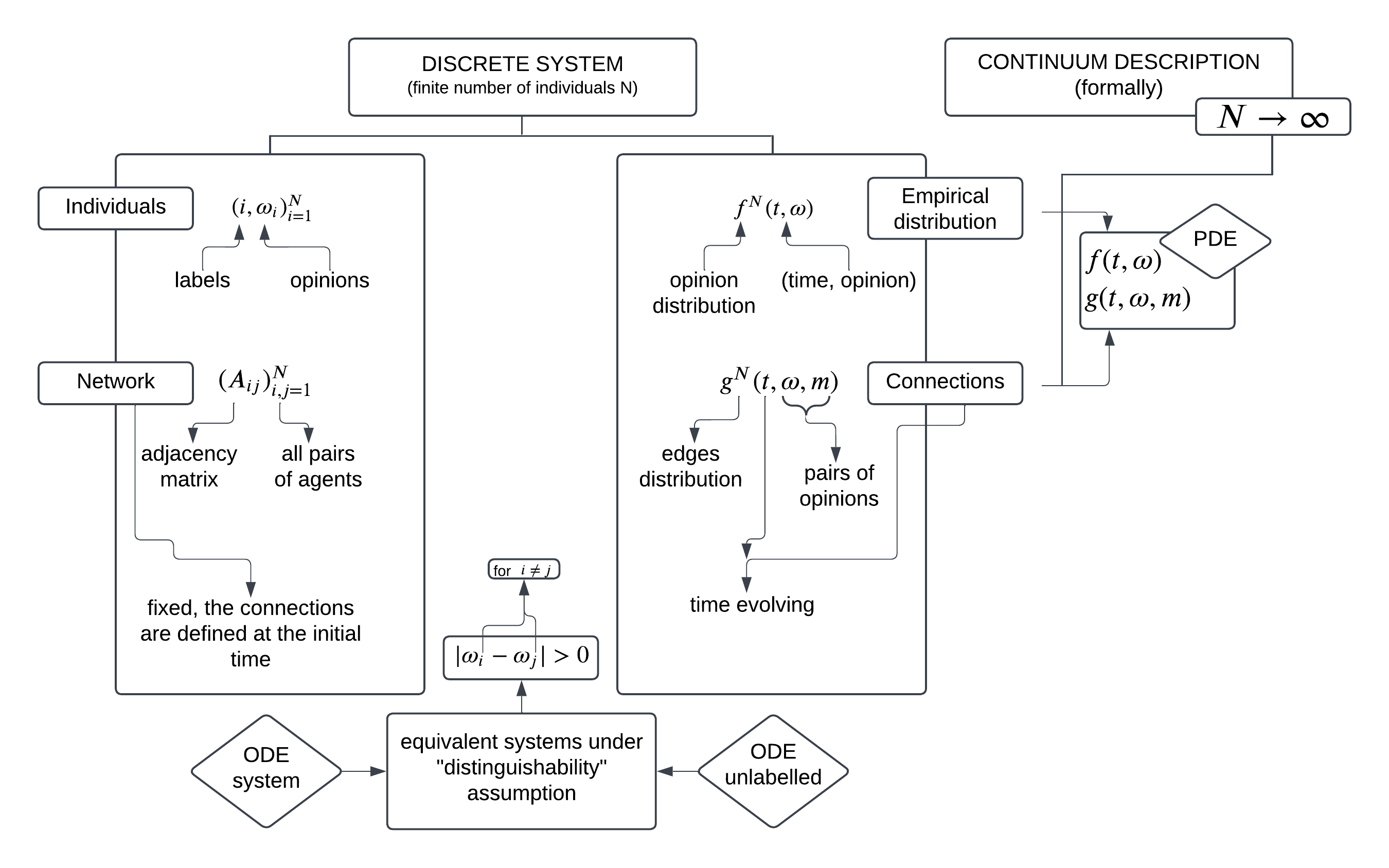}
    \caption{
    \textbf{Diagram of the modeling framework.} The discrete system (left) is presented from two perspectives: one describes the evolution of individual opinions \((\omega_i)_{i=1}^N\) over time, where interactions occur on a network represented by the adjacency matrix \((A_{i,j})_{i,j=1}^N\) (see Sec.~\ref{sec:iba}). The other perspective considers the distributions of opinions, \(f^N\), and edges, \(g^N\), as described in Sec.~\ref{sec:network description}. Both perspectives yield equivalent sets of ODEs for short times, provided the "distinguishability" assumption \ref{as:initial_time} holds (corresponding to Lemma~\ref{lem:ODE_approx}). From \(f^N\) and \(g^N\), the continuum model is formally derived by taking the limit as \(N \to \infty\) (Sec.~\ref{sec:limit}).
    }
    \label{fig:diagram}
\end{figure}

Our goal here is to establish a description for the dynamics of the population that is as simple as possible so that concrete computations and comparison with the discrete dynamics can be carried out. Moreover, this description should be amenable to adding changes to the network so that it can capture dynamical networks, as well as adding noise, and other types of forces (see Section~\ref{sec:extensions_model}). The approach presented here is based on introducing a density distribution of edges, and it is adapted from Degond et al.~\cite{degond2016continuum} The overall picture of this framework can be found in Fig. \ref{fig:diagram}.

%
%
%
%
%
%
%
%

\subsection{Network description: node distribution and edge distribution}
\label{sec:network description}
The goal of this Section is to find a continuum model that approximates the dynamics of~\eqref{eq:ODE_eb} in the limit of a large  number of individuals $N\to \infty$. In particular, we would like to determine the time-evolution for $f=f(t,\omega)$, the distribution of agents with opinion $\omega$ at time $t$.

Large-particle limits are the central theme of mean-field limits. However, model~\eqref{eq:ODE_eb}, rather than having a mean-field force, has a network-mediated force (to have a mean-field force the equation should have the sum  over all the individuals and not just over the individuals to which one is connected). As a consequence, we cannot apply directly the methodology of mean-field limits. The main bottleneck to apply this methodology is that the network is characterized using the labels of the agents, but one would need that the interactions between the agents can be described with independence of labels. For this reason, we give next a characterization of the network using the opinions of the agents instead of their labels. In this sense, every agent $i$ will be characterized only by its opinion $\omega_i$ instead of by its opinion \textit{and} its label.\\
Specifically, instead of having a node (individual) $i$ with value $\omega_i$, we consider a node $\omega_i$; and instead of having an edge $\{i,j\}$, we consider an edge  $\{\omega_i, \omega_j\}$.  The function
\begin{align}\label{eq:bijection_labels}
\phi\,:\,(i\mapsto \omega_i(t))_{i=1}^N \quad \mbox{is a bijection for} \quad t\in[0,T]
\end{align}
between $\{1,\hdots,N\}$ and $\{\omega_1(t),\hdots, \omega_N(t)\}_{t\in [0,T]}$,
as long as $|\omega_i(t)-\omega_j(t)|>r>0$ for $i\neq j$ and $t\in [0,T]$ and for some $r$. Thanks to this, we are able to characterize uniquely all the individuals (labels) by just keeping track of the opinions.

For this reason, we make the following assumption:
\begin{assumption}[Distinguishability assumption] \label{as:initial_time}
    At $t=0$, for all $i,j\in \{1,\hdots, N\}$, $i\neq j$ it holds for some $r>0$ that
    $$|\omega_i(0)- \omega_j(0)|>2r>0.$$
\end{assumption}
This assumption extends at least for short times, as explained in the following Lemma:
\begin{lemma} \label{lem:assumption}
Suppose that Assumption~\ref{as:initial_time} holds, then
for some $T>0$, and for all times $t\in [ 0,T]$ and all $i,j\in \{1,\hdots, N\}$, $i\neq j$ it holds that
    $$|\omega_i(t)- \omega_j(t)|>r.$$
In particular,~\eqref{eq:bijection_labels} holds and, $0<r=r(N)$ is such that
\begin{equation}
\label{eq:def_r}
r\leq\inf_{t\in[0,T]}\min_{i\neq j}|\omega_i(t)-\omega_j(t)|.
\end{equation}
\end{lemma}
Lemma~\ref{lem:assumption} is a direct consequence of the continuity of solutions over time.

Since we want to have a continuum description for the interactions in the network, we need to obtain a description for the distribution of the nodes and the distribution of the edges.  First, we define the node distribution as the  empirical distribution of opinions given by
\begin{equation}\label{eq:def_fN}
f^N(t,\omega):=\frac{1}{N}\sum_{i=1}^N \delta_{\omega_i(t)}(\omega),
\end{equation}
where $\delta_{\omega_i}$ denotes the Dirac delta distribution centered at $\omega_i$.
Our goal is to determine the time-evolution for $f^N$ as $N\to\infty$. Denote by $B_{r/2}(\omega)$ the open ball of radius $r/2$ centered at $\omega$ and consider the integral:
$$I(t,\omega_0):=\int_{B_{r/2}(\omega_0)}\xd f^N(t,\omega), \qquad \omega_0 \in \Omega,\, \quad t\in [0,T].$$
Under the notations and the assumptions of Lemma~\ref{lem:assumption} we have that $I(t,\omega_0)=1/N$ if at time $t$ there is an agent with opinion $\omega_0$ and 0 otherwise (by Lemma~\ref{lem:assumption} we know that if there is an agent with opinion $\omega_0$, then this agent is unique, and it is at distance at least $r$ from other opinions). Therefore, in this setup, $f^N$ is able to distinguish between individuals. For this reason, we say that $f^N$ not only captures the distribution of the opinions, but also the nodes of the network.

To have a full description of the network,
we need one more ingredient: since we describe the network through the opinions,  we also need an evolution for the \emph{description} of the network. The \emph{network} itself is fixed, but the opinions are not, since the opinions evolve over time, the description of the network also evolves over time. This is the price that we have to pay for getting rid of the labels in~\eqref{eq:ODE_eb}. With this end in mind, we introduce the distribution of edges $g^N$ given by
\begin{align} \label{eq:def_gN}
g^N(t,\omega,m):=\frac{1}{2\kappa}\sum_{i\sim j} \delta_{\omega_i(t)}(\omega) \delta_{\omega_j(t)}(m),
\end{align}
where $\kappa$ is the total number of connections (or edges) at initial time (but since the network is fixed, it remains constant). We divide by 2 because in the sum we are counting the pairs twice. Notice that $g^N$ is symmetric in its variables.

Proceeding analogously as before, for $t\in [0,T]$ and $(\omega_0,m_0) \in \Omega^2$ we have that the integral $$\Gamma(t,\omega_0,m_0):=\int_{B_{r/2}(\omega_0)}\int_{B_{r/2}(m_0)}\xd g^N(t,\omega,m)$$
satisfies $\Gamma(t,\omega_0,m_0)=\kappa^{-1}$ if at time $t$ there exist two individuals with opinions $\omega_0$ and $m_0$ (and if this pair exists, it is unique) that are connected via the network; and $\Gamma(t,\omega_0,m_0)=0$ otherwise. Therefore, from $g^N$ we are able to recover the network. \\
Moreover, the integral
\begin{align*}
\int_{B_{r/2}(\omega_0)}\int_{\Omega} \xd g^N(t,\omega,m)
\end{align*}
is 0 if there is no individual with opinion $\omega_0$ and, if it is different from 0, then it means that there is exactly one individual with that opinion. In this case, the integral corresponds to the total number of connections of $\omega_0$ divided by $2\kappa$. Thus, from this information, one can also recover the information on the nodes, and hence, reconstruct the node density $f^N$, but, in general, we will have that
\begin{equation}
    \label{eq:g is not f}
    \int g^N(t,\omega, dm) \neq c f^N(t,\omega), \qquad \mbox{ for any constant } c>0\,,
\end{equation}
the equality being the exception: it occurs when every individual is connected with all the others. In this case $c= N/(2\kappa)$). To illustrate that fact, we have the following remark.

\begin{remark}[Example]
\label{rem:example irrecoverability f from g}
A simple example (illustrated in Figure~\ref{fig:diag network}) of this consists of a system with 3 individuals with opinions $\omega_1 < \omega_2 < \omega_3$, with connections $1 \leftrightarrow 2$ and $2 \leftrightarrow 3$, so that $\kappa=2$. We have then
$$f^N(\omega)= \frac{1
}{3}\left(\delta_{\omega_1}(\omega)+\delta_{\omega_2}(\omega)+\delta_{\omega_3}(\omega)\right)$$
and
$$
    g^N(\omega,m) = \frac{1}{4}\bigg( \delta_{\omega_1}(\omega)\delta_{\omega_2}(m)+\delta_{\omega_2}(\omega)\delta_{\omega_1}(m)  +\delta_{\omega_2}(\omega)\delta_{\omega_3}(m)+\delta_{\omega_3}(\omega)\delta_{\omega_2}(m)\bigg)\,.
$$
Additionally
$$\int g^N(\omega,dm) = \frac{1}{4}\left( \delta_{\omega_1}(\omega)+2\delta_{\omega_2}(\omega)+\delta_{\omega_3}(\omega)\right),$$
which is not proportional to $f^N$.

This example also highlights the importance of distinguishable opinions: if two individuals have the same opinion, we can no longer differentiate between them or accurately describe the network.
\end{remark}

\begin{figure}
	\centering

	\begin{tikzpicture}[scale=3, node/.style={circle, draw, minimum size=7mm, inner sep=0pt}]

		\newcommand\wa{0.25}
		\newcommand\wb{0.6}
		\newcommand\wc{0.8}
		\newcommand\gap{1.0}

		\node[node] (w1) at (-1.2-\gap, 0.15) {$1$};
		\node[node] (w2) at (-0.8-\gap, 0.75) {$2$};
		\node[node] (w3) at (-0.4-\gap, 0.15) {$3$};

		\draw[Stealth-Stealth] (w1) -- (w2);
		\draw[Stealth-Stealth] (w2) -- (w3);

		\draw[-Stealth] (-0.1, 0) -- (1, 0) node[right] {$\omega$};
		\draw[-Stealth] (0, -0.1) -- (0, 1) node[left=1pt] {$m$};

		\draw (\wa, -0.02) node[below] {$\omega_1$} -- (\wa, 0);
		\draw (\wb, -0.02) node[below] {$\omega_2$} -- (\wb, 0);
		\draw (\wc, -0.02) node[below] {$\omega_3$} -- (\wc, 0);
		\draw[dotted,darkgray] (\wa, -0.0) -- (\wa, 1);
		\draw[dotted,darkgray] (\wb, -0.0) -- (\wb, 1);
		\draw[dotted,darkgray] (\wc, -0.0) -- (\wc, 1);

		\draw (-0.02, \wa) node[left] {$\omega_1$} -- (0, \wa);
		\draw (-0.02, \wb) node[left] {$\omega_2$} -- (0, \wb);
		\draw (-0.02, \wc) node[left] {$\omega_3$} -- (0, \wc);
		\draw[dotted,darkgray] (0, \wa) -- (1, \wa);
		\draw[dotted,darkgray] (0, \wb) -- (1, \wb);
		\draw[dotted,darkgray] (0, \wc) -- (1, \wc);

		\node[circle,red] at (\wa, \wb) {$\bullet$}; 
		\node[circle,red] at (\wb, \wa) {$\bullet$}; 
		\node[circle,red] at (\wc, \wb) {$\bullet$}; 
		\node[circle,red] at (\wb, \wc) {$\bullet$}; 

		\draw [gray, decorate, decoration = {brace,raise=5pt}] (\wb,\wa) -- (\wa,\wb) node[black, anchor=center,pos=0.5,below left = 16pt and -7pt,rotate=-45] {$1 \leftrightarrow 2$};
		\draw [gray, decorate, decoration = {brace,raise=5pt}] (\wb,\wc) -- (\wc,\wb) node[black, anchor=center,pos=0.5,above right = 16pt and -7pt,rotate=-45]{$2 \leftrightarrow 3$};

	 	\node at (-0.8-\gap, 1.2) {Network of connections};
		\node at (0.5, 1.2) {Representation of $g^N(\omega, m)$};


	\end{tikzpicture}

  \caption{Illustration of the example in Remark \ref{rem:example irrecoverability f from g}.
    On the right, braces correspond to an edge in the network, and a \textcolor{red}{$\bullet$} represents a Dirac delta. These come in pairs, as the network is not directed.
	}
	\label{fig:diag network}
\end{figure}

\subsubsection{Smoothing of the edge distribution}
In this Section, all the integrals are integrals on the domain $\Omega$.

In the sequel, we will need to evaluate the distribution $g^N$ at some specific opinion $\omega$, which is not possible since we have delta-distributions. For this reason, we introduce a smooth approximation $g_r^N$ of $g^N$:
\begin{equation}
 g_r^N(\omega,m):=\frac{r^2}{2\kappa}\sum_{i\sim j} \delta^{(r)}_{\omega_i(t)}(\omega)\delta^{(r)}_{\omega_j(t)}(m)   = \frac{r^2}{2\kappa} \sum_{i=1}^N\sum_{j\in \mathcal{I}_i}\delta^{(r)}_{\omega_i(t)}(\omega)\delta^{(r)}_{\omega_j(t)}(m), \label{eq:definition_gN}
\end{equation}
 where $\delta^{(r)}_{\omega_0}$ is a mollifier, which is a smooth approximation to the $\delta$-distribution centered at $\omega_0\in \Omega$, i.e. $\delta^{(r)}_{\omega_0}\to \delta_{\omega_0}$ as $r\to 0$ in distribution. In particular, we will assume that
\begin{align} \label{eq:assumption_delta}
&\delta^{(r)}_{\omega_0}(\omega) \mbox{ compactly supported in } [\omega_0-r/2, \omega_0+r/2] \nonumber\\
& \int \delta^{(r)}_{\omega_0}=1, \mbox{ and } \delta^{(r)}_{\omega_0}(\omega_0) = \frac{1}{r}.
\end{align}
 Notice that $r^{-2}g^N_r$ is a probability distribution on $\Omega\times \Omega$ and that
 $$\frac{1}{r^2} g^N_r(\omega,m)\to g^N(\omega,m) \quad \mbox{as} \quad r\to 0$$
 in distribution.

\subsection{An individual-based model using the distributions}
\label{sec:IBM_nolabels}

 We can rewrite our original ODE system \eqref{eq:ODE_eb} as follows:

\begin{lemma} \label{lem:ODE_approx}
Under Assumptions~\ref{as:initial_time}, for $t<T$, the system of equations \eqref{eq:ODE_eb} is equivalent to
\begin{align} \label{eq:mean-field-force}
    \frac{d\omega_i(t)}{dt}= a[g^N_r](\omega_i(t)) + B^{(N)}_r(\omega_i(t)),
\end{align}
where $B^{(N)}_r\to 0$ as $r\to 0$ for fixed $N$ ($B^{(N)}_r$ is given in \eqref{eq:B_definitions}),
and we defined the operators
\begin{align}
  a[g](\omega):=& \int \eta[g](\omega,m) \mathbf{D}(\omega-m) \xd m,\nonumber
        \\
        \eta[g](\omega,m) :=& \frac{g(\omega,m)}{\int g(\omega,m') \xd m'},\label{def:eta}
\end{align}
with $\eta[g](\omega, \cdot)=0$ if $g(\omega,\cdot) =0$ a.e..
\end{lemma}

\begin{proof}[Proof of Lemma~\ref{lem:ODE_approx}]
  Thanks to introducing the approximation of the $\delta$-distributions, we can evaluate $g_r^N$ at a given point $\omega_i(t)$ using~\eqref{eq:assumption_delta}:
 \begin{align}
    g_r^N(\omega_i(t),m) &= \frac{r^2}{2\kappa}\sum_{j\in \mathcal{I}_i} \delta^{(r)}_{\omega_i(t)}(\omega_i(t)) \, \delta^{(r)}_{\omega_j(t)}(m)= \frac{r}{2\kappa}\sum_{j\in \mathcal{I}_i} \delta^{(r)}_{\omega_j(t)}(m)\,.\nonumber
    \intertext{Integrating with respect to $m$ and using the fact that the $\delta_{\omega_j(t)}^{(r)}$ have disjoint support (Lemma~\ref{lem:assumption}), we obtain}
\label{eq:g omega_i integrated}
    \int g_r^N(\omega_i, m) \xd m &= \frac{r}{2\kappa} \sum_{j\in\mathcal{I}_i} \int \delta_{\omega_j(t)}(m) \xd m\,,
    \intertext{ and since all integrals in the sum equal one, we have}
\label{eq:expr card Ii continuous}
    \#\mathcal{I}_i&:= \frac{2\kappa}{r}\int g_r^N(\omega_i(t),m)\,\xd m.
\end{align}
We can then rewrite \eqref{eq:ODE_eb}:
\begin{align*}
  \frac{\xd \omega_i}{\xd t} &= \frac{1}{\# \mathcal{I}_i} \sum_{j\in \mathcal{I}_i} \mathbf D(\omega_i - \omega_j)
     = \frac{1}{\# \mathcal{I}_i }\sum_{j\in \mathcal{I}_i}\int \mathbf D( \omega_i-m) \,\delta_{\omega_j}(m)\, \xd m \\
    & =\frac{1}{\# \mathcal{I}_i } \sum_{j\in \mathcal{I}_i} \int \mathbf D(\omega_i-m) \left( \delta_{\omega_j}(m) -\delta^{(r)}_{\omega_j}(m)\right) \xd m + \frac{1}{\# \mathcal{I}_i }\sum_{j\in \mathcal{I}_i} \int \mathbf D(\omega_i-m) \delta^{(r)}_{\omega_j}(m) \xd m  =: {B}_r^{(N)} + {C}_r^{(N)},
\end{align*}
with
\begin{align}
 B_r^{(N)}&:= \frac{1}{\# \mathcal{I}_i } \sum_{j\in \mathcal{I}_i} \int \mathbf D(\omega_i-m) \left( \delta_{\omega_j}(m) -\delta^{(r)}_{\omega_j}(m)\right)\xd m \nonumber\\
 &= \frac{r}{2\kappa \int g^N_r(\omega_i(t),m) \xd m} \sum_{j\in \mathcal{I}_i} \int \mathbf D(\omega_i-m) \left( \delta_{\omega_j}(m) -\delta^{(r)}_{\omega_j}(m)\right)\xd m \label{eq:B_definitions} \to 0 \mbox{ as } r\to 0\,.
 \end{align}
Using~\eqref{eq:g omega_i integrated} and \eqref{eq:expr card Ii continuous} we also have
\begin{align*}
   C_r^{(N)}&:=\frac{1}{\# \mathcal{I}_i }\sum_{j\in \mathcal{I}_i} \int\mathbf D(\omega_i-m) \delta^{(r)}_{\omega_j}(m) \xd m
                =\frac{1}{\# \mathcal{I}_i}\frac{2\kappa}{r} \int \mathbf D(\omega_i-m)  g_r^N(\omega_i, m) \xd m\\
                  &=\int \mathbf D(\omega_i-m) \frac{g_r^N(\omega_i, m)}{\int g_r^N(\omega_i, m') \xd m'} \xd m= \int \mathbf D(\omega_i-m) \eta[g_r^N](\omega, m)\xd m.
\end{align*}
For $t<T$, it thus holds that
\begin{align}
  \frac{d \omega_i}{dt} =  B^{(N)}_r + \int \mathbf \eta[g_r^N](\omega, m) \mathbf{D}( \omega_i-m) \xd m,
\end{align}
which yields the result.
\end{proof}

\subsection{\texorpdfstring{Large-particle limit $N\to \infty$}{Large-particle limit N to infinity}}
\label{sec:limit}

Our goal is to compute formally the limit $N\to \infty$ for the model \eqref{eq:mean-field-force}.
\begin{remark}[Equivalence of the systems]
 To have that the dynamics of \eqref{eq:ODE_eb} and \eqref{eq:mean-field-force} match, we must consider $r=r(N)$ as given in Lemma~\ref{lem:assumption}. Then we know that the dynamics for both systems coincide for $t<T$, $T=T(r, N)$. However, forcibly, it must hold that $r(N)\to 0$ as $N\to \infty$. Therefore, in the limit, we cannot assume that Lemma~\ref{lem:ODE_approx} holds. However, for a fixed given $r$, the system of ODEs~\eqref{eq:mean-field-force} is well-defined and we can take the limit $N\to\infty$.
 \end{remark}

\begin{assumption} \label{as:limit}
    Let $f^N$ and $g_r^N$ be defined in~\eqref{eq:def_fN} and~\eqref{eq:definition_gN}, respectively. Suppose that both functions have limits as $N\to \infty$ (and $r(N)\to 0$), denoted by $f$ and $g$  --- suppose that the convergence is strong enough, and the limiting functions are regular enough (so that we can take the limit and do integration by parts in \eqref{eq:aux_limit} and in \eqref{eq:aux_limit2} below). Under these assumptions, we take $r(N)$ small enough so that
\begin{align}
\label{eq:assumption_remainder_term}
\sum_{i=1}^N |B^{(N)}_{r(N)}| \to 0 \quad \mbox{ as } \quad N\to \infty.
\end{align}
(This is possible by choosing an approximation of the $\delta$-distribution $\delta_r$ that converges to $\delta$ fast enough.)
\end{assumption}

\begin{remark}
   Proving rigorously the limit of $f^N$ is a central objective in the study of mean-field limits. Various techniques have been developed for this purpose, such as coupling approaches (see~\cite{carmona2016lectures,sznitman1991topics}). For a comprehensive treatment, we refer the reader to the monographs~\cite{chaintron2021propagation,chaintron2022propagation}.

In contrast, establishing the limit of $g_r^N$ is significantly more challenging. This is the primary difficulty on making this methodology rigorous. It is reasonable to conjecture that the limit of $g_r^N$ requires certain scaling assumptions related to the network size $N$.

Since the computations that follow are formal, we will assume that limits can be exchanged with derivatives and integrals. In this context, when we require the convergence to be "sufficiently strong", we mean that we assume that the convergence of the functions is such that the formal manipulations performed in the subsequent calculations are justified.
\end{remark}

\begin{proposition}[Large-particle limit for \eqref{eq:mean-field-force}]
\label{prop:mean-field-limit}
    Consider $f^N$ and $g_r^N$ under Assumption~\ref{as:limit}. We denote by $f=f(t,\omega)$ and $g=g(t,\omega,m)$ their respective limits. Then, these functions satisfy the following system of equations:
\begin{subequations}
  \label{eq:kinetic system}
\begin{align}
  \partial_t f(\omega) &+ \partial_\omega \left(a[g](\omega) f(\omega)\right) = 0\,, & \omega \in \Omega\,,
    \label{eq:mean-field opinion}\\
\partial_t g(\omega, m)&+ \partial_\omega \left(a[g](\omega) g(\omega,m)\right)
    + \partial_m \left(a[g](m) g(\omega,m)\right)=0\,, & (\omega, m) \in \Omega^2\,,
    \label{eq:mean-field connectivity}
    \end{align}
\end{subequations}
where $a[g]$ is given in Lemma~\ref{lem:ODE_approx}.\\ The initial conditions $f(t=0)$ and $g(t=0)$ are given as the limits of $f^N(t=0)$ and $g^N(t=0)$, respectively.\\
To ensure the preservation of the total mass, we impose zero flux boundary conditions at the boundary of the domain $\partial \Omega$:
\begin{subequations}
\label{eq:kinetic system_boundary}
    \begin{align}
    a[g](\omega)f(\omega) &= 0\,, & \omega & \in \partial\Omega \label{eq:mean-field bc opinion} \\
    a[g](\omega)g(\omega,m)&= 0\,, & (\omega,m) & \in \partial\Omega^2 \label{eq:mean-field bc network}\, .
    \end{align}
\end{subequations}
\end{proposition}

The proof can be found at the end of this Section --- we comment first on this result.
Equation \eqref{eq:mean-field opinion} gives the evolution for the distribution of opinions $f=f(t,\omega)$. Notice that the equation is in a conservative form, which guarantees that the total mass of individuals is preserved. It can be interpreted as a continuity equation with ``velocity'' $a[g](\omega)$, i.e. the ``velocity'' of opinions is dictated by $a[g]$, which describes interactions through the network. The function $\eta(\omega,\cdot)$ in~\eqref{def:eta} is a probability distribution that gives the distribution of opinions connected to opinion $\omega$. Notice that $g$ is symmetric but $\eta$ is not.

Since the network is represented via the opinions and the opinions evolve over time, the distribution $g=g(\omega,m)$ also evolves over time. Hence, we get equation \eqref{eq:mean-field connectivity} which is again an equation in conservative form (which guarantees that the total density of links is constant over time). It has two terms, one for each one of the variables $\omega$ and $m$. These terms are symmetric so $g$ is also symmetric, i.e. $g(\omega,m)=g(m,\omega)$. We observe that the equation for $g$ is closed (i.e., it is independent of $f$). Moreover, if we define the connectivity of an opinion $\omega$ as
\begin{align}\label{eq:def_connectivity}
h(\omega):= \int g(\omega, m)\, \xd m,
\end{align}
we observe that $h$ and $f$ satisfy the same equation (one obtains the equation for $h$ by integrating the equation of $g$ against $m$). However, $h$ and $f$ do not have the same initial data, so they are not the same. In particular, this means that from $g$ we cannot recover $f$ (recall \eqref{eq:g is not f} at the discrete level).

\medskip
Further, notice that in the mean-field limit as it is carried out here, we lose the information of the relation between $\kappa$, which we recall, it represents the total number of connections, and $N$. Moreover, because of the 0-homogeneity of $\eta$ in $g$ (i.e. $\eta[\lambda g]=\eta[g]$ for any scalar $\lambda \neq 0$),
it holds the following:
\begin{proposition}[Scaling invariance]
    If $(f, g)$ is a solution of system \eqref{eq:kinetic system} with initial data $(f_0, g_0)$,
    then, $(\lambda f, \mu g)$ is a solution with initial data $(\lambda f_0, \mu g_0)$, for any $\lambda, \mu > 0$.
\end{proposition}
Even though we do not exploit this property in the current paper, it is a relevant feature of the differential equation that can be used in future mathematical analyses.

The proof is a direct check.
In particular, by taking $\lambda = 1$ as required by $f$ being a probability distribution, this means that multiplying the initial datum $g(t=0)$ by a constant will leave $f$ unchanged.
In other words, the density $f$ remains invariant under scaling of the initial data $g(t=0)$. However, in our case we fix the total integral of $g(t=0)$ to be 1 to have a probability density.

\begin{proof}[Prop. \ref{prop:mean-field-limit}]
Let us consider a smooth test function $\varphi$ compactly supported in the interior of $\Omega$ and compute the time derivative of
\begin{align}
\frac{\xd}{\xd t}\Bigl(\int \varphi \xd f^N \Bigr)=&\frac{\xd}{\xd t} \left(\frac{1}{N}\sum_{i=1}^N\varphi(\omega_i(t))\right)=\frac{1}{N}\sum_{i=1}^N \varphi'(\omega_i(t))\frac{\xd \omega_i}{\xd t}  \nonumber\\
=&\frac{1}{N}\sum_{i=1}^N\varphi'(\omega_i(t))\,  a[g^N_r](\omega_i(t)) + \frac{1}{N}\sum_{i=1}^N \varphi'(\omega_i(t)) B^{(N)}_r(\omega_i(t)) \nonumber \\
=& \int \varphi'(\omega) \, a[g^N_r](\omega) df^N(t,\omega) + \frac{1}{N}\sum_{i=1}^N \varphi'(\omega_i(t)) B^{(N)}_r(\omega_i(t))\,. \label{eq:aux_limit}
\end{align}
Letting $N\rightarrow\infty$, using \eqref{eq:assumption_remainder_term} and integrating by parts, since we assume that $f$ and $g$ are regular enough, we obtain
$$
    \frac{\xd}{\xd t}\Bigl( \int\varphi \xd f \Bigr)=\int \varphi'(\omega) a[g](\omega)\, f(\omega)\xd \omega    =-\int \varphi(\omega)\partial_\omega[ a[g](\omega) f(\omega)]\,\xd \omega \, ,
$$
which corresponds to the equation of $f$ \eqref{eq:mean-field opinion} in weak form.

We need a second equation which describes the evolution of $g$.
Let $\varphi=\varphi(\omega,m)$, $\varphi\in C^\infty_c(\Omega\times\Omega)$ be a test function with compact support in the interior of the domain and consider
\begin{align}
\int\varphi(\omega,m)\xd g^N_r(\omega,m)=\frac{r^2}{2\kappa}\sum_{i\sim j}\varphi(\omega_i,\omega_j). \label{eq:aux_limit2}
\end{align}
Now, computing the time derivative and integrating by parts, we obtain (proceeding as before) the equation \eqref{eq:mean-field connectivity} for $g$.
\end{proof}

\subsection{Case of consensus dynamics: properties captured from the discrete system}
\label{sec:consensus continuum}
In this Section, we perform an analysis parallel to the one done for the discrete model (Section~\ref{sec:consensus_discrete}) with the aim of checking if the essential structure of the discrete dynamics \eqref{eq:ODE_eb} is captured by the continuum model~\eqref{eq:mean-field opinion}--\eqref{eq:mean-field connectivity} when we have consensus dynamics (Assumption~\ref{ass:debate operator}).

First, we have the analogous conserved quantity of Lemma~\ref{lem:conserved_quantity_discrete}:
\begin{lemma}\label{lem:conserved_quantity_continuum}
It holds that
 \begin{align} \label{eq:conserved_quantity_continuum}
   \int\int \omega g(t,\omega, m) \xd\omega \xd m =\int\int \omega g(0,\omega, m) \xd\omega \xd m.
 \end{align}
\end{lemma}
The proof is a direct check (by time-differentiating). Notice that, indeed, expression \eqref{eq:conserved_quantity_continuum} is the analog to the one in the discrete system given in \eqref{eq:dt}.

Second, we can define the analog of the discrete potential $V$ \eqref{eq:potential_V} for the continuum dynamics as:
$$
\tilde V(t)=\iint W(\omega- m) g(t, \omega, m)\xd \omega\xd m.
$$

Third, we obtain a similar ``Lyapunov condition'' \eqref{eq:Lyapunov_condition_micro} as in the discrete setting for the continuum potential $\tilde V$:

\begin{lemma}
 It holds that
  \begin{align}
    \frac{\xd}{\xd t}\tilde V(t)
    \leq - 2\int h(\omega) a^2[g](\omega)\, d\omega \leq 0 \label{eq:estimate_second_moment_g}
  \end{align}
 where $h$ is defined in \eqref{eq:def_connectivity}.
\end{lemma}

\begin{proof}
We compute the time derivative of $\tilde V$:
  \begin{align*}
    \frac{\xd}{\xd t}&\left( \int\int W(\omega-m) g(t,\omega, m) \xd \omega \xd m\right) \\
    = &- \int\int W(\omega-m) \partial_\omega (a[g](\omega) g(\omega,m)) \xd\omega \xd m- \int\int W(\omega-m) \partial_m (a[g](m) g(\omega,m)) \xd\omega \xd m\\
     = &- 2 \int\int W(\omega-m) \partial_\omega (a[g](\omega) g(\omega, m)) \xd\omega \xd m
     =  - 2 \int\int \mathbf{D}(\omega-m) a[g](\omega) g(\omega,m) d\omega dm\\
     = & - 2 \int h (\omega) a^2[g](\omega) d\omega,
  \end{align*}
  where in the second equality we used that $W(\omega-m) = W(m-\omega)$ and $g(\omega,m)=g(m,\omega)$; and in the third one we used integration by parts, the zero-flux boundary condition \eqref{eq:mean-field bc network}, and that $\mathbf{D}(\omega-m) = - \partial_\omega W(\omega-m)$. The last inequality is obtained by multiplying and dividing the integrand by $h(\omega)$ and using that $\eta(\omega,m)=g(\omega,m)/h(\omega)$ and the definitions.
\end{proof}

 In analogy with the discrete system, we will suppose that
\begin{equation} \label{eq:assume_limit_fg}
 f(t,\omega) \to \delta_{\oinfcont}(\omega)  \text{ and }
 g(t, \omega,m) \to \delta_{\oinfcont}(\omega)\delta_{\oinfcont}(m)
\end{equation}
 as $t\to \infty$ for some constant $\oinfcont$, which is the final consensus opinion.
 \begin{remark}
 Proving the convergence to the equilibrium for the continuum system is more challenging than in the discrete one. Notice that it is even not clear how to check that the limiting distributions are stationary solutions of the continuum equations \eqref{eq:mean-field opinion}--\eqref{eq:mean-field connectivity} since one has to divide by $\int g(\omega,m) \xd m$.
 \end{remark}
Assuming the limits \eqref{eq:assume_limit_fg},  it holds that
$$
\int \int \omega g(\omega,m) \to \oinfcont,
$$
and by applying Lemma~\ref{lem:conserved_quantity_continuum}, the limiting value $\oinfcont$ is given by:
\begin{equation}
\oinfcont =  \int\int \omega g(t=0, \omega, m) \xd\omega \xd m.
\label{eq:limit omega continous}
\end{equation}
Finally, we notice that the limiting value $\oinfcont$ is the analog to the discrete limiting value $\oinfmicro$ in \eqref{eq:limit_omega}.

In conclusion, some structural properties of the discrete system investigated in Section~\ref{sec:consensus_discrete} are captured by the continuum dynamics.

\subsection{Conclusions from the basic continuum modeling framework}
From the proposed modeling framework, we conclude that the continuum system \eqref{eq:mean-field opinion}-\eqref{eq:mean-field connectivity} is well-suited for both numerical simulations and mathematical analysis. Numerical simulations will be presented in Section~\ref{sec:numerics}, while in the previous Section, we provided a basic mathematical analysis in the context of consensus dynamics, verifying that certain properties of the discrete system are preserved in the continuum model.\\
In this regard, we have successfully achieved our objective of developing an application-oriented framework for modeling continuum equations in opinion dynamics. However, an important limitation remains: we have not rigorously derived these continuum equations from the discrete system. Instead, we constructed the model using inspiration from mean-field approximations.\\
Notably, Lemma~\ref{lem:assumption} loses validity as $N \to \infty$, which can result in a loss of crucial network information in the continuum limit, as demonstrated in Remark~\ref{rem:validity lemma}

\begin{figure}
  \begin{tikzpicture}
    \draw (0, 0) node[inner sep=0] {\includegraphics[width=0.45\textwidth]{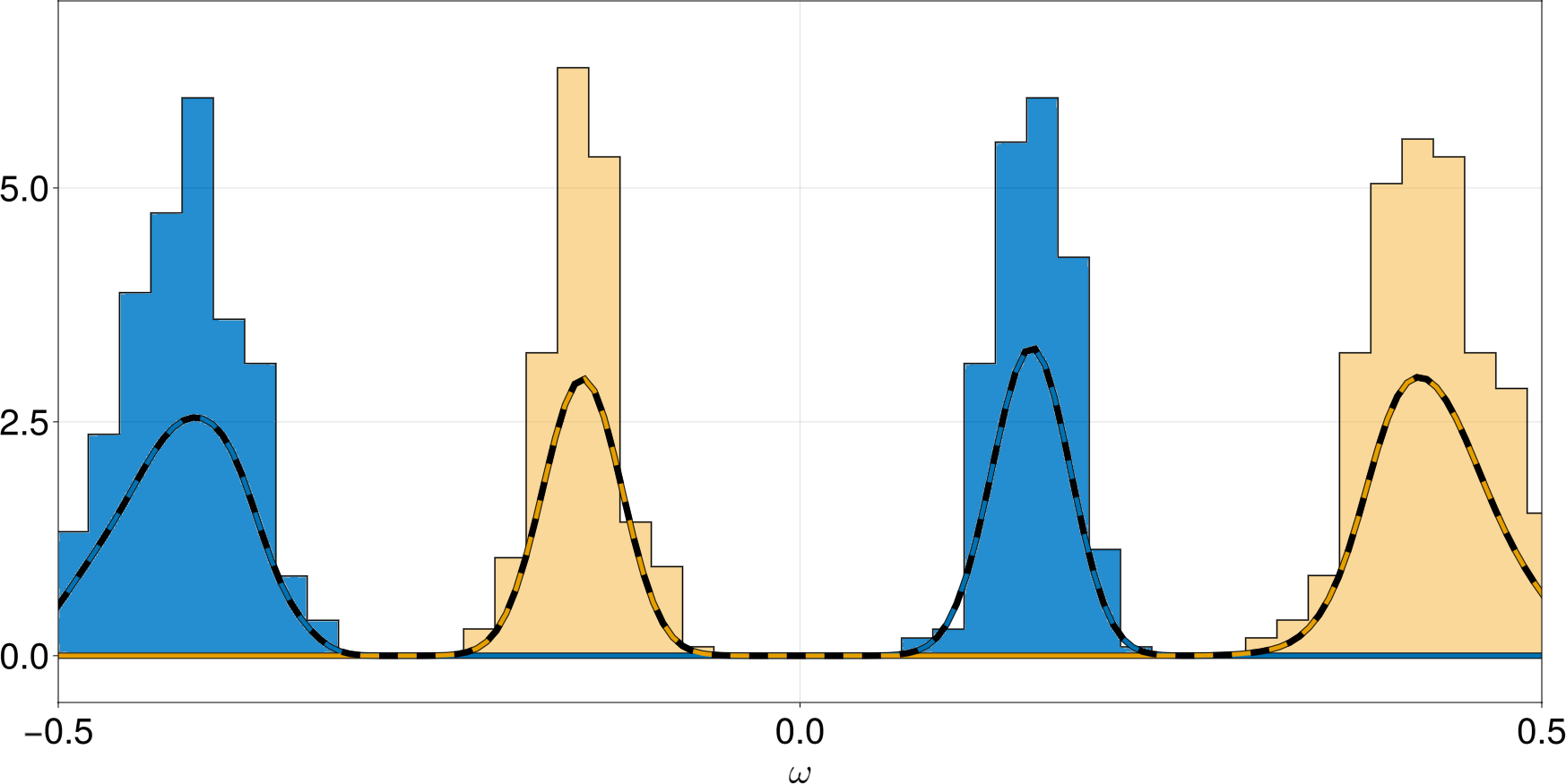}};
    \draw (-2.0, 1.4) node {$t=0$};
    \draw (-2.4, 0.7) node {$B^\textrm{out}$};
    \draw (-2.5, -0.7) node {$\longrightarrow$};
    \draw (-1, 0.7) node {$Y^\textrm{in}$};
    \draw (-0.7, -0.7) node {$\longrightarrow$};
    \draw (0.7, 0.7) node {$B^\textrm{in}$};
    \draw (0.9, -0.7) node {$\longleftarrow$};
    \draw (3.0, 0.7) node {$Y^\textrm{out}$};
    \draw (2.7, -0.7) node {$\longleftarrow$};
  \end{tikzpicture}
  \begin{tikzpicture}
    \draw (0, 0) node[inner sep=0] {\includegraphics[width=0.45\textwidth]{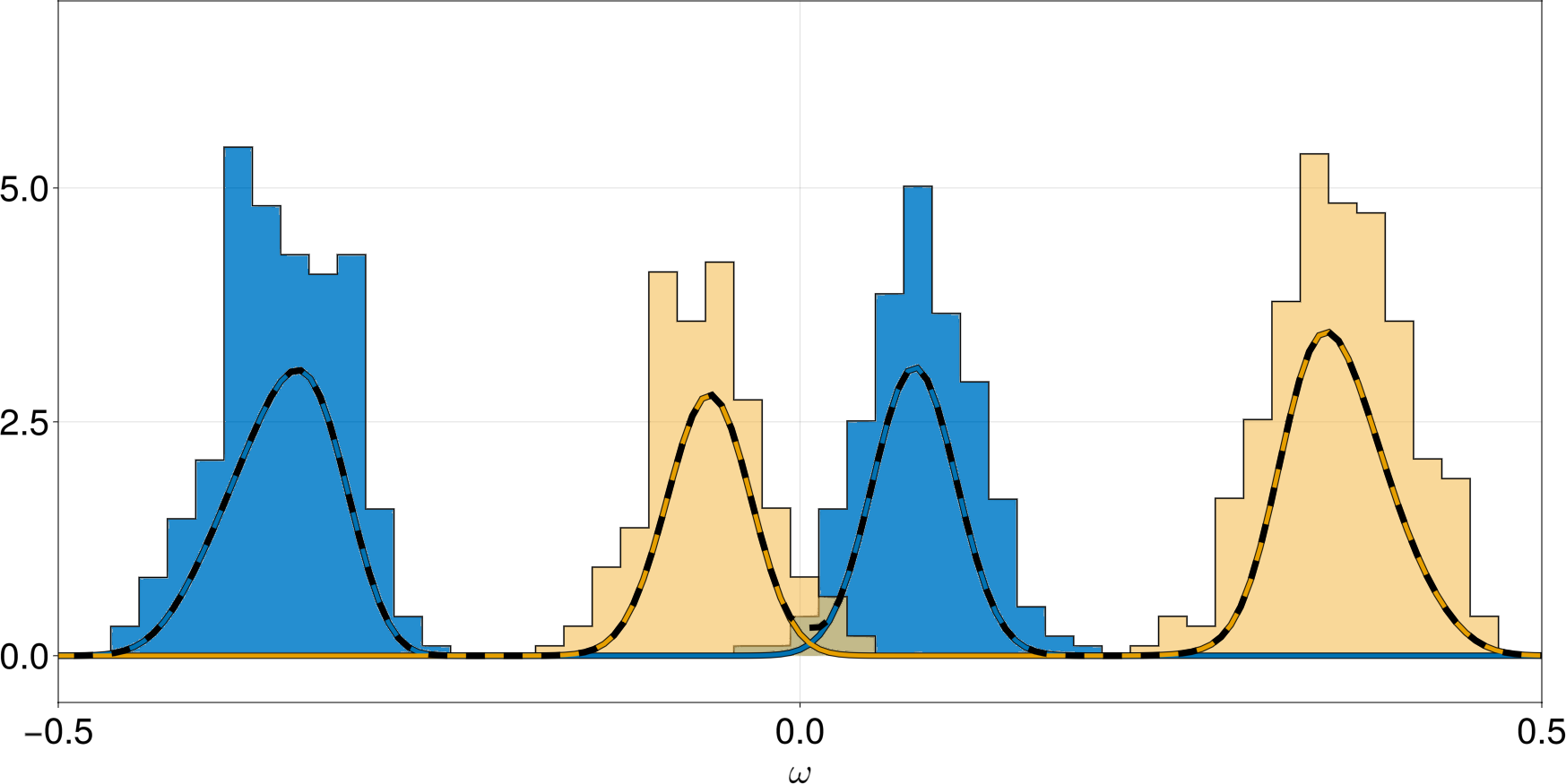}};
    \draw (-2.5, 1.4) node {$t=0.3$};
  \end{tikzpicture}
  \begin{tikzpicture}
    \draw (0, 0) node[inner sep=0] {\includegraphics[width=0.45\textwidth]{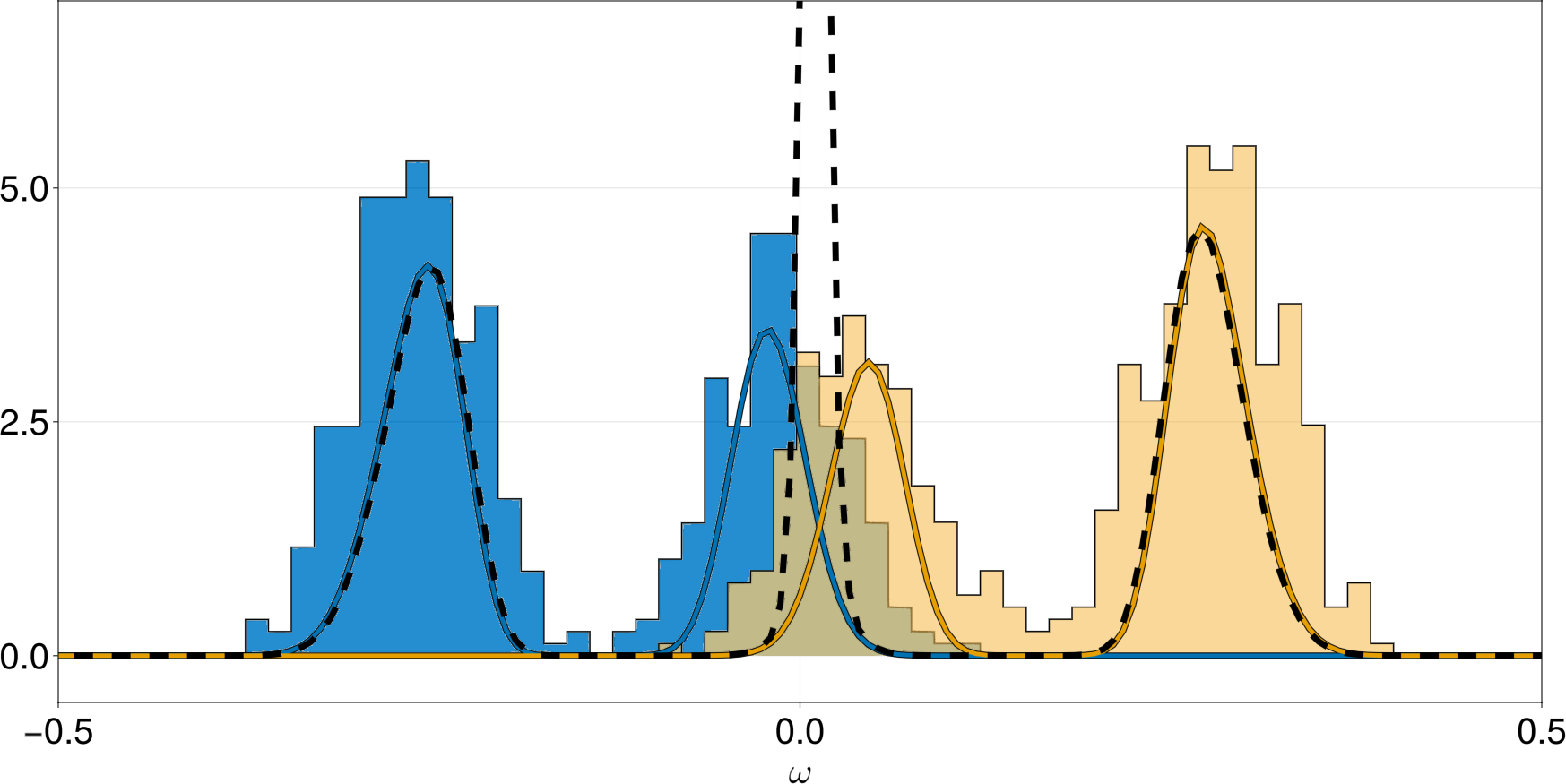}};
    \draw (-2.5, 1.4) node {$t=0.9$};
    \draw (1, 1.4) node {$\longleftarrow$ spike};
  \end{tikzpicture}
  \begin{tikzpicture}
    \draw (0, 0) node[inner sep=0] {\includegraphics[width=0.45\textwidth]{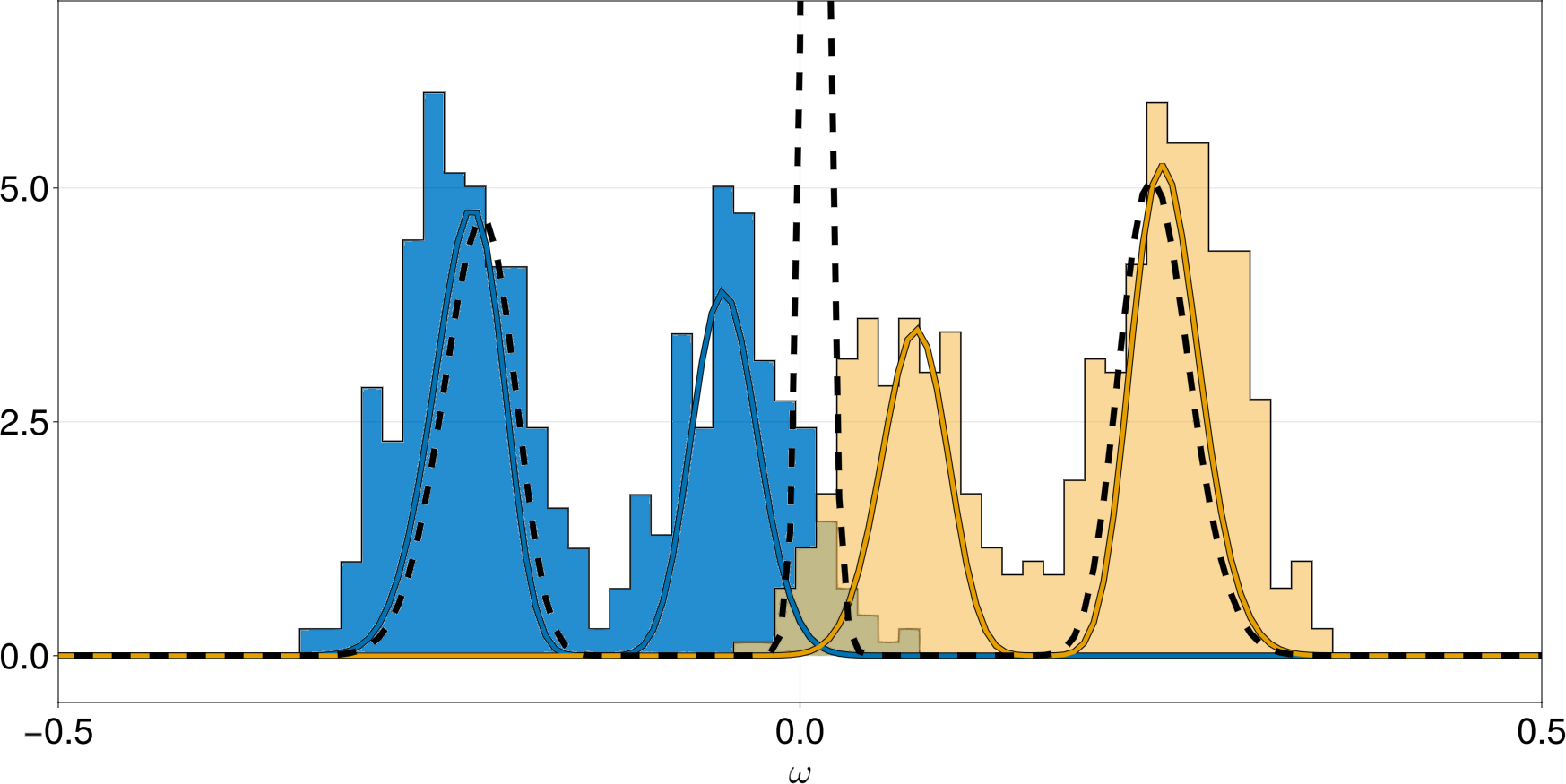}};
    \draw (-2.5, 1.4) node {$t=1.2$};
    \draw (0.7, 0.7) node {$Y^\textrm{in}$};
    \draw (0.9, -0.7) node {$\longrightarrow$};
    \draw (-0.55, 0.7) node {$B^\textrm{in}$};
    \draw (-0.5, -0.7) node {$\longleftarrow$};
  \end{tikzpicture}
  \caption{Comparison between the discrete dynamics (histograms), the continuum model without group labeling (dash line) and the continuum model with group labeling (solid lines).
    The histogram represents the evolution of the discrete solution, while the dashed line is the solution to the continuous model without group labeling. Here, $N=1000$. See Remark~\ref{rem:validity lemma} for details.
  }
   \label{fig:snapshots}
\end{figure}

\begin{remark}\label{rem:validity lemma}
Here we consider a simple example where the basic continuum (without group labeling) system fails to capture the dynamics of the discrete model.
See Figure~\ref{fig:snapshots} for an illustration.

\emph{Setup:} The population is split in two groups, Blue (B) and Yellow (Y).
The two groups are poorly connected to each other (the mixing parameter is $\mu=1/20$, in other words, each individual has on average only 1 in 20 of its connections going to the other group). The initial distribution in each group is not homogeneous, and consist of two subgroups denoted \emph{in} and \emph{out} respectively.

\emph{Dynamics}: Within a given group, as the agents are well connected,
the distribution homogenizes quickly. Initially, there is a good agreement between the discrete and continuous opinion distributions, but starting at around $t=0.3$,
there is an overlap between subgroups $B^\textrm{in}$ and $Y^\textrm{in}$, which cross and separated later in time in the microscopic model.
This is not captured correctly at the continuous level without group labeling: the subgroups merge and stop drifting around the origin, leading to a stationary spike in the distribution (see the third frame in Figure~\ref{fig:snapshots}).
On the contrary, the continuum model \emph{with group labeling} presented in the following, captures the discrete dynamics well (see the solid colored lines).
\end{remark}

To address this issue, we propose remedies in the following Section. These remedies involve augmenting the system with additional information to retain the network structure. However, this comes at the cost of increased complexity, potentially making the system less tractable for both analysis and simulation.

\section{Continuum modeling framework: refinements}
\label{sec:refinements}

\paragraph{The regime of validity of Lemma \ref{lem:ODE_approx}.}

In this article, we have relied on the validity of Lemma~\ref{lem:ODE_approx}, which establishes the equivalence between the original discrete system \eqref{eq:ODE_eb} and the model \eqref{eq:mean-field-force}. For finite \( N \), this equivalence holds, provided that \(\omega_i(t) \neq \omega_j(t)\) for all \(i, j = 1, \ldots, N\) over some finite time interval. In other words, the equivalence is preserved as long as individuals can be distinguished by their opinions. An example of a set-up that would not preserve this property is given in Fig. \ref{fig:snapshots}.

To overcome the challenge of indistinguishable individuals (i.e., when \(\omega_i(t) = \omega_j(t)\) for some \(i \neq j\)), additional information can be introduced into the system. In the original paper~\cite{degond2017continuum}, which inspired our derivation, the authors considered a system of particles defined by their positions and orientations. In that framework, particles interact only when they are spatially close, and spatial repulsion prevents them from overlapping, allowing them to remain distinguishable by using the information of their positions.

In our case, where spatial variables are absent, we can introduce an alternative form of distinguishing information. For example, instead of tracking only an agent's opinion \(\omega_i(t)\) at time \(t\), we could also track their initial opinion \(\omega_i(0)\). This way, if the opinion of two agents overlap at some time \(t\), they can still be distinguished by their distinct initial opinions.

Consequently, the continuum dynamics would give the evolution of the opinion distribution \( f = f(t,\omega,\omega_0) \), which gives the density of agents with opinion \(\omega\) at time \(t\) who initially held opinion \(\omega_0\) at time \(t=0\). Similarly, the edge distribution would be given by \( g = g(t,\omega,\omega_0,m,m_0) \), representing the proportion of edges between agents with opinions \(\omega\) and \(m\) at time \(t\) who initially had opinions \(\omega_0\) and \(m_0\), respectively (see remark \ref{rem:incorporating_initial_data} below for more details).

While this approach preserves the ability to distinguish individuals, it comes at the cost of increasing the dimensionality of the system, thereby making computations significantly more costly.

\paragraph{Incorporating community structure.}
\label{sec:with community structure}

To simplify the problem and facilitate numerical testing, we propose a different approach: incorporating community information. Specifically, we assume the social network follows an LFR-type structure, where individuals are divided into $n>0$ communities. Each agent \(i\) is represented by the pair \((\omega_i, p_i)\), where \( \omega_i \) is the agent's opinion and \( p_i \) is a community label, which is time-independent. The opinion dynamics remain the same as in equation \eqref{eq:ODE_eb}, while the community label evolves trivially as
\[
\frac{dp_i}{dt} = 0.
\]

With this modification, we define the empirical opinion density and the edge distribution as:
\begin{align}
&f^N(t, \omega, p):= \frac{1}{N}\sum_{i=1}^N \delta_{(\omega_i(t),p_i)}(\omega, p),\\
&g^N_r(t, \omega, p, m, q):= \frac{r^2}{2\kappa}\sum_{i\sim j}\delta_{(\omega_i(t), p_i)}(\omega,p)\delta_{(\omega_j(t),p_j)}(m, q),
\end{align}
where \( p \) and \( q \) are community labels.

To distinguish individuals with the same opinion but different group memberships, we define the distance between two individuals \((\omega, p)\) and \((m, q)\) as:
\[
\|(\omega, p) - (m, q)\| = |\omega - m| + |p - q|.
\]
We impose the following assumption on the initial configuration:
\begin{assumption}
\label{as:assumption_group}
    There exists $r=r(N)>0$ such that
    $$\inf_{i\neq j}\|(\omega_i(0), p_i)-(\omega_j(0), p_j)\| > 2r(N)>0.$$
\end{assumption}
Under this assumption, we can extend the results of Lemma~\ref{lem:ODE_approx}. Specifically, we have:
\begin{lemma} Under assumption \ref{as:assumption_group} there exists $T>0$, such that for all $t\in [0,T]$, the discrete system of equations \eqref{eq:ODE_eb} is equivalent to:
\begin{align}\label{eq:discrete_groups}
\frac{d\omega_i(t)}{dt}= a[g^N_r](\omega_i(t), p_i) + B_r^{(N)}(\omega_i(t), p_i),
\end{align}
where $B_r^{(N)}(\omega_i(t), p_i)\to 0$ as $r\to 0$ for fixed $N$, and
\begin{align}
& a[g^N](\omega_i, p_i):= \sum^n_{q=1}\int \mathbf{D}(\omega_i-m) \eta^N(t,\omega_i,p_i,m,q) dm,\\
&\eta^N(t,\omega,p, m,q) := \frac{g^N(t,\omega,p,m,q)}{\sum_{q=1}^n\int g^N(t,\omega,p,m,q)},
\end{align}
with $\eta[g](\omega,p,\cdot,\cdot)=0$ if $g(\omega,p,\cdot,\cdot)=0$ a.e..
\end{lemma}

Following analogous assumptions and procedure, we derive the continuum equations for \eqref{eq:discrete_groups}:
\begin{align} \label{eq:modified_macro1}
    &\partial_t f(t,\omega,p) + \partial_\omega\left( a[g](\omega,p) f(t,\omega,p) \right)=0,\\
    & \partial_t g(\omega,p,m,q)+\partial_\omega\left(a[g](\omega,p)g(\omega,p,m,q)\right) + \partial_m\left( a[g](m,q) g(\omega,p,m,q)\right)=0,\\
    & a[g](\omega,p) f(t,\omega,p) =0 \qquad\qquad\omega \in \partial \Omega,\\
    &a[g](\omega,p)g(\omega,p,m,q)=0 \qquad  (\omega,m) \in \partial \Omega^2. \label{eq:modified_macro4}
\end{align}

Equations~\eqref{eq:modified_macro1}--\eqref{eq:modified_macro4} are structurally equivalent to the original system \eqref{eq:mean-field opinion}-\eqref{eq:mean-field connectivity}, but the functions \( f \) and \( g \) now depend on the additional variables \( p \) and \((p, q)\). This seemingly minor modification significantly increases the computational complexity.

For example, in our numerical simulations with \(n \) communities, we must compute $n$ functions for the opinion density \(f(t, \omega, p)\) (one for each group) and $n^2$ functions for the edge distribution \(g(t, \omega, p, m, q)\), corresponding to all pairs \(p, q \in \{1, \hdots, n\}\). In contrast, the original system only required solving two equations, highlighting the increase in complexity from \(2\) to \(n + n^2 \) equations.

\begin{remark} \label{rem:incorporating_initial_data}
If we had chosen to track the initial opinions \(\omega_0\) and \(m_0\) instead of community labels, the resulting continuum system would still correspond to \eqref{eq:modified_macro1}-\eqref{eq:modified_macro4}, but having $\omega_0$ instead of $p$ and $m_0$ instead of $q$. But notice that \(\omega_0\) and \(m_0\) take continuous values in \(\Omega\): this dramatically increases the computational cost, as opposed to handling discrete community labels \(p\) and \(q\).
\end{remark}

\paragraph{Conclusion: tailoring modeling choices.}

We conclude this Section by emphasizing that the continuum modeling framework must often be tailored to the specific system under consideration. Depending on the application, additional information such as community structure, initial data, or other distinguishing variables may need to be incorporated to accurately capture the network's dynamics. These modifications, while necessary, increase the complexity of both the mathematical model and its numerical implementation.

\section{Numerical comparison: consensus dynamics between three communities}
\label{sec:numerics}
\newcommand{\boldW}{\boldsymbol{\omega}}
This Section is dedicated to numerical experiments for both the discrete model of Section~\ref{sec:iba} and the continuous system of Section~\ref{sec:mfl}.
We start by giving some details about the discretization, both in time and in opinion space, followed by considerations about the choice of initial data,
and conclude by a comparison of both models for specific choices for the Lancichinetti–Fortunato–Radicchi (LFR) network~\cite{LFR}, explained below.

\subsection {Discretization}
\paragraph{The discrete model.}
We consider \eqref{eq:ODE_eb} with a fixed number of agents $N$ on a given graph, described by the adjacency matrix $A$.
For the time integration, we use an explicit Euler scheme: let $\Delta t > 0$ be the (fixed) time-step, and let $t_n = n \; \Delta t$,
$n \in \N_{\ge 0}$ be the corresponding discrete times. We can then define the time-discrete opinions as
\[ \omega_i^n := \omega_i(t = t_n)\, \quad 1 \le i \le N\,.\]
Equation \eqref{eq:ODE_eb} is then integrated in time using
\begin{equation}
  \label{eq:micro:discretized}
  \omega_i^{n+1}
= \omega_i^n + \Delta t \frac{1}{\# \mathcal{I}_i} \sum_{j\in \mathcal{I}_i} \mathbf{D}(\omega_i^n - \omega_j^n)\,,
\end{equation}
where $\Delta t$ must be chosen small enough to guarantee $\omega_i^{n+1} \in \Omega$, regardless of Assumptions~\ref{ass:debate operator}.
  More precisely, if $\Delta t \le \| \partial_z \mathbf{D}\|_{L^\infty}^{-1}$, one can prove that $\omega_i^{n+1}$ in the convex hull of $\{\omega_i^n\} \cup \{ \omega_j^n : j \in \mathcal{I}_j\}$,
  and in particular $\omega_i^{n+1} \in \Omega$. For details, see e.g. Appendix~\ref{sec:app:stepsize}.

\paragraph{The continuum system.}
\label{sec:numerics mfl}
\newcommand{\Ncont}{N_\mathrm{cont}}
We now turn to the system \eqref{eq:kinetic system}, which we solve using an adaptation of the classical local Lax–Friedrichs (LLF) finite volume scheme~\cite{leveque2002finite}.
The Lax–Friedrichs method was initially proposed to approximate solutions of hyperbolic, linear conservation laws, later generalized to the nonlinear case, that is to equations of the form
\begin{align}
\label{eq:general_form_finite_volume}
\partial_t q + \partial_\omega(b(q(\omega))) = 0\,, \qquad \omega \in \Omega\,,
\end{align}
there the flux function $b$ is a function of \emph{pointwise} values of the unknown $q$, which makes the equation \emph{local}.
Although it is only first order accurate and suffers from high numerical diffusion, it is conservative and is relatively easy to extend.

System~\eqref{eq:kinetic system} falls out of the scope of the LLF method for at least two reasons: first, it is a system of two equations,
one set on $\Omega$ and coupled, through the factor $a[g]$, to the other, which is set on $\Omega \times \Omega$.
Second, not only equations~\eqref{eq:kinetic system} are nonlinear, but the flux $\omega \mapsto a[g](\omega) f(\omega)$ is not a function of $f(\omega)$ alone. It depends \emph{non locally} on $g$, and is \emph{space varying}.
The integral factor $a[g]$ thus introduces space dependence, coupling and non locality.

The LLF scheme also covers non-constant flux functions, so the first point is not a problem.
We can deal with the second point by neglecting the dependency of $a[g]$ on $g$, i.e. by first fixing $a$, then solving \eqref{eq:mean-field opinion} and \eqref{eq:mean-field connectivity} independently at each time-step with the LLF scheme and finally updating $a$.
This way, we are back in the linear setting.

More precisely, we fix the discretization with $\Ncont > 0$ and define $\Delta \omega = |\Omega|/\Ncont$. For $0 \le i \le \Ncont-1$, define the finite volumes $I_i = [-1 + i \Delta \omega, -1 + (i+1) \Delta \omega]$ and integrate~\eqref{eq:mean-field opinion} on $I_i$ to get
\[ \frac{\xd}{\xd t} \int_{I_i} f + \int_{\partial I_i} a f = 0\,. \]
For some measurable set $B$, we define the average $f|_B = \int_B f / |B|$, so that we can now take the piecewise constant approximation of $f$, $f_i := f|_{I_i}$. Now, $a_i$ is defined similarly, by averaging $a[g]$.
This yields
\begin{equation}
  \label{eq:llf f}
\frac{\xd}{\xd t} \int_{I_i} f_i + \widehat{a f}_{i+1/2} - \widehat{a f}_{i-1/2}   = 0\,,
\end{equation}
where $\widehat{a f}_{i \pm 1/2}$ is the so-called numerical flux, which takes the place of the (ill-defined) boundary terms.
Likewise for equation~\eqref{eq:mean-field connectivity}, we get the following integrated equation for $g_{i,j}$:
\begin{equation}
  \label{eq:llf g}
  \frac{\xd}{\xd t} \int_{I_i \times I_j} g_{i,j} + \widehat{a g}_{i+1/2,j} - \widehat{a g}_{i-1/2,j}  + \widehat{a g}_{i,j+1/2} - \widehat{a g}_{i,j-1/2}  = 0\,.
\end{equation}
We take the local Lax–Friedrichs flux, which in our case is given by
\[ \widehat{a f}_{i+1/2} = \frac{1}{2} \left(a_if_i + a_{i+1} f_{i+1} -  (f_{i+1} - f_i) \max(|a_i|,|a_{i+1}|)\right)\,.\]
\emph{Discretization of the velocity $a$:} here, we have neglected the dependence of $a$ on the pointwise value of $f$ on $I_i$, which is reasonable for large values of $\Ncont$.
We use simple Riemann sums to compute $\eta_i$ and $a_i$, i.e. we first compute
\[\eta_j :=
  \begin{cases}
    g_{jk} / (\Delta \omega \sum_{k} g_{jk}) & \text{ if } \Delta \omega \sum_{k} g_{jk} < \epsilon
  \\
  0 & \text{ otherwise},
  \end{cases}
\]
where $\epsilon$ is some small cutoff parameter (in practice, $\epsilon = 10^{-10}$).
The discrete velocity is then defined as \[a_i := \Delta \omega \sum_j \eta_j D_{ij}\,,\] where $D_{ij} = \mathbf{D}(\midint(I_i) - \midint(I_j))$ and $\midint(I)$ is the middle point of interval $I$.

We can use a similar approach for $g$, where we need to distinguish between the two directions (note the index for $a$):
\begin{align*}
\widehat{a\,g}_{i+\frac{1}{2},j} &= \frac{1}{2} \left(a_i g_{i,j} + a_{i+1}g_{i+1,j} -  (g_{i+1,j} - g_{i,j}) \max(|a_i|, |a_{i+1}|)\right)\,,
\\
\widehat{a\,g}_{i,j+\frac{1}{2}} &= \frac{1}{2} \left(a_j g_{i,j} + a_{j+1}g_{i,j+1} -  (g_{i,j+1} - g_{i,j}) \max(|a_j|, |a_{j+1}|)\right)\,,
\end{align*}
the remaining fluxes $\widehat{a f}_{i - \frac{1}{2}}$, $\widehat{a\,g}_{i-\frac{1}{2},j}$ and $\widehat{a\,g}_{i,j-\frac{1}{2}}$ are defined likewise.

The boundary conditions \eqref{eq:mean-field bc opinion}-\eqref{eq:mean-field bc network} are enforced by setting
\begin{gather*}
  \widehat{a f}_{-1/2} = \widehat{a f}_{\Ncont - 1/2} = 0\,,
  \\
  \widehat{a\,g}_{-1/2, j} = \widehat{a\,g}_{\Ncont - 1/2, j} = \widehat{a\,g}_{i, -1/2} = \widehat{a\,g}_{j, \Ncont - 1/2} = 0\,, \quad 0 \le i,j \le \Ncont - 1\,.
\end{gather*}

We then discretize \eqref{eq:llf f} and \eqref{eq:llf g} in time with time-step $\Delta t$ to get
\begin{align}
    f_i^{n+1} - f_{i}^n &= - \frac{\Delta t}{\Delta \omega}\left(\widehat{a f}_{i+1/2}^n - \widehat{a f}_{i-1/2}^n\right)\,,
    \label{eq:f discrete evolution}
    \\
    g_{i,j}^{n+1} - g_{i,j}^n &= - \frac{\Delta t}{\Delta \omega}\big(\widehat{a\,g}_{i+1/2,j}^n - \widehat{a\,g}_{i-1/2,j}^n + \widehat{a\,g}_{i,j+1/2}^n - \widehat{a\,g}_{i,j-1/2}^n\big) \,,
    \label{eq:g discrete evolution}
\end{align}
where the exponent $n$ indicates the approximation at time $t^n = n \Delta t$.

\begin{remark}[Properties of the numerical scheme]
  From the point of view of numerical analysis, \eqref{eq:kinetic system} is complex in several ways: it is a system of hyperbolic equations, set on different spatial domains, one of which is two dimensional. The velocity $a$ is also non local, depending on both space and time, and also lacks any kind of regularity. As such, the system does not fit ---to our knowledge--- any existing analytical framework. A rigorous study of it falls outside of the scope of the current paper.

  Nonetheless, this simple numerical scheme has a number of good properties.
  \begin{itemize}
    \item It as a finite volume scheme, it is conservative by construction.
    \item It preserves the symmetry of $g$, i.e. if $g_{i,j}^n = g_{j,i}^n$ then $g_{i,j}^{n+1} = g_{j,i}^{n+1}$, this follows from the symmetry of the right-hand side of \eqref{eq:g discrete evolution}.
      Swapping $i$ and $j$ in the third term gives $\widehat{a\,g}_{j,i+1/2}^n$ which equals the first term $\widehat{a\,g}_{i+1/2,j}^n$, as can be easily checked. The remaining terms are dealt with similarly.
      Note that this is independent of how $a_i$ is defined.
    \item Under the CFL-type condition $\Delta t < \frac{\Delta \omega}{2 ||\mathbf{D}||_\infty}$, it preserves the positivity of $f$ and $g$. Consider the case of $f$, and assume that $f_i^n > 0$ for all $i$.
      Starting from \eqref{eq:f discrete evolution} and the definition of the numerical fluxes we have
      \begin{align*}
        0 < f_i^{n+1} &\Leftrightarrow \widehat{a f}_{i+1/2}^n - \widehat{a f}_{i-1/2}^n < \frac{\Delta \omega}{\Delta t} f_i^n \\
                      & \Leftrightarrow \frac{1}{2} \Big(\underbrace{(a_{i+1} - \max(|a_{i+1}|, |a_i|) f_{i+1}^n}_{\le 0} + \underbrace{(a_{i-1} - \max(|a_i|, |a_{i-1}|) f_{i-1}^n}_{\le 0}  \\
                      & \qquad + \underbrace{(\max(|a_{i+1}|,|a_i|) + \max(|a_i|,|a_{i-1}|))}_{\le \max_i |a_i|} f_i^n \Big) < \frac{\Delta \omega}{\Delta t} f_i^n\,.
      \end{align*}
      From its definition, $a_i$ is a convex combination of the $D_{ij}$, so that $|a_i| \le \| \mathbf{D} \|_\infty$
      so that $\Delta t < \frac{\Delta \omega}{\| \mathbf{D} \|_\infty}$ is a sufficient condition. The factor $\frac{1}{2}$ in the condition above is necessary when considering $g$, to account for fluxes in both directions.
  \end{itemize}
\end{remark}


\begin{remark}
This numerical scheme is extended in a straightforward manner when we consider the continuum equations containing the group labels \eqref{eq:modified_macro1}-\eqref{eq:modified_macro4}.
\end{remark}

\subsection{Choice of the network and initial opinion distribution}
\label{sec:numerics initial data}

  The specific choice of initial opinion distribution and more importantly the choice
  of the network will depend on the real system being modeled. There are many possibilities to do that, and it is out of the scope of the present
  paper to review them. We will focus on interactions between communities using the LFR network.

  Since our goal here is to assess the continuous system's ability to approximate the discrete model,
  we construct the initial data of the continuous system from that of the discrete as follows:
  \begin{enumerate}
    \item fix the number of agents $N$ for the discrete equation \eqref{eq:ODE_eb} and the discretization size $\Ncont$ for the continuous system~\eqref{eq:kinetic system};
    \item pick an initial opinion distribution $f_0$ and a network $G$ of size $N$;
    \item construct the initial opinions $(\omega_i^0)_{1 \le i \le N}$ by sampling $f_0$;
    \item compute the initial discrete distribution $(f_i^0)_{0 \le i \le \Ncont-1}$ by averaging $f_0$ over $I_i$;
    \item compute the initial discrete distribution $(g_{i,j}^0)_{0 \le i,j \le \Ncont-1}$ by kernel density estimation (KDE) of $g^N(t=0)$ as defined in \eqref{eq:def_gN} and averaging over $I_i \times I_j$. The kernel bandwidth is taken equal to the (one-dimensional) Sheather-Jones bandwidth~\cite{Sheater-91} corresponding to $\omega_i^0$.
      We use the KDE as a black box to get a continuous approximation of the sum of Dirac deltas. We refer the interested reader to Chen's~\cite{chen_tutorial_kde_2017} introduction to the method.
\end{enumerate}

For the choice of network, we are particularly interested in how the presence of communities affect the agreement between both models.
For this, we consider graphs and opinion distributions initially displaying some amount of so-called homophily.
That is: clusters in the network structure corresponding to localized bumps in the opinion distribution, and vice versa.
By changing how isolated the clusters are within the network and how localized the opinions are within these clusters, one may construct different initial conditions.

\begin{figure}
    \includegraphics[width=\linewidth]{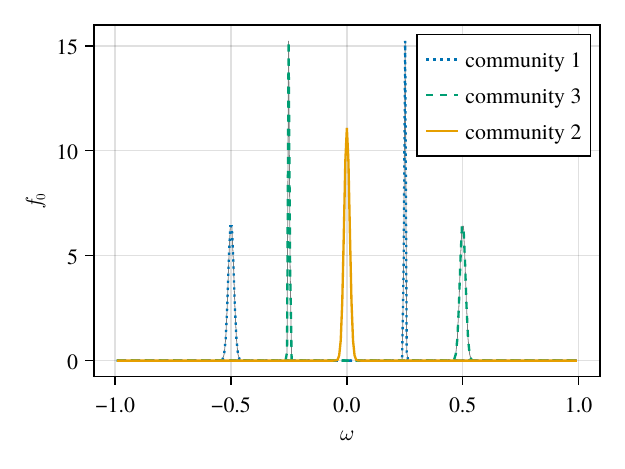}
    \caption{Initial opinion distribution, assuming three communities of identical size.  Each color corresponds to a community.}
    \label{fig:f0 var}
\end{figure}

In practice, we use Lancichinetti–Fortunato–Radicchi's algorithm~\cite{LFR} to generate graphs, where each node (or individual) $i$ is assigned to a community $c(i) \in \{1,2,3\}$.
Within each of these communities, the initial distribution of opinion is similar to that of Figure~\ref{fig:snapshots}. More precisely, the initial distribution is a weighted sum of normal distributions, centered on $-1/2$ and $1/4$ for community 1, $0$ for community 2, and $-1/4$ and $1/2$ for community 3, with the \emph{outer} normal components (those centered on $-1/2$ and $1/2$) having a higher weight and variance than the \emph{inner} ones. This is arguably an artificial setup, but it helps to highlight the impact of introducing group labeling in the continuous model; this will become clearer below.
Communities do not have equal sizes in general, and since the opinion distribution is weighted by community size, one cannot expect the overall distribution to be even as a function.
The overall initial distribution $f_0$ is then a weighted sum of five truncated normal distributions, weighted by community size. The corresponding initial opinion distributions are shown in Figure~\ref{fig:f0 var}.

\begin{figure}
    \centering
    \begin{subfigure}[b]{0.3\linewidth}
        \includegraphics[width=\linewidth]{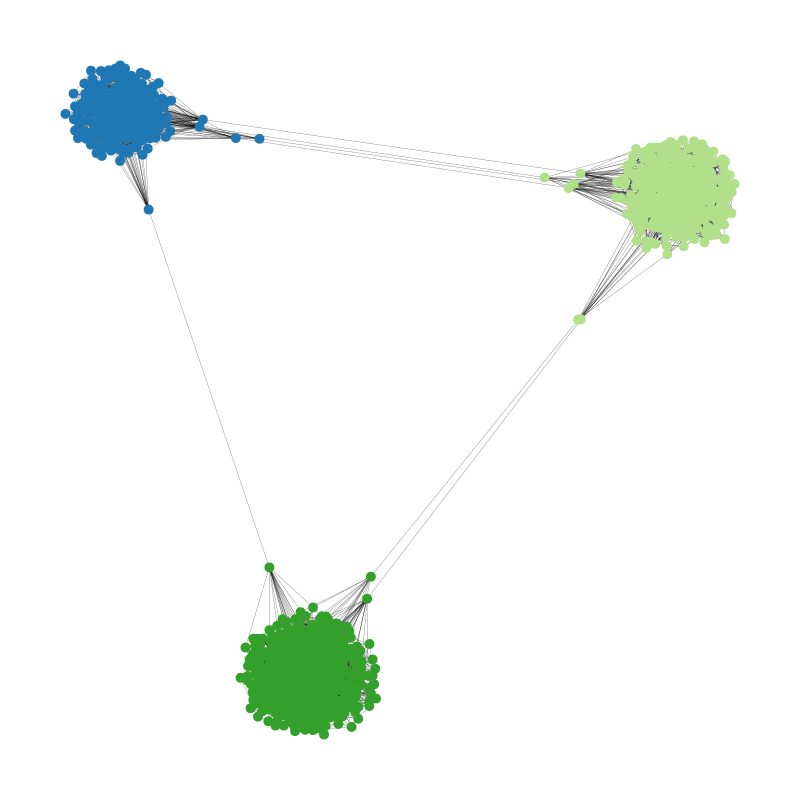}
        \caption{$\mu = \num{e-3}$}
        \label{fig:lfr graphs 0.001}
    \end{subfigure}
    \begin{subfigure}[b]{0.3\linewidth}
        \includegraphics[width=\linewidth]{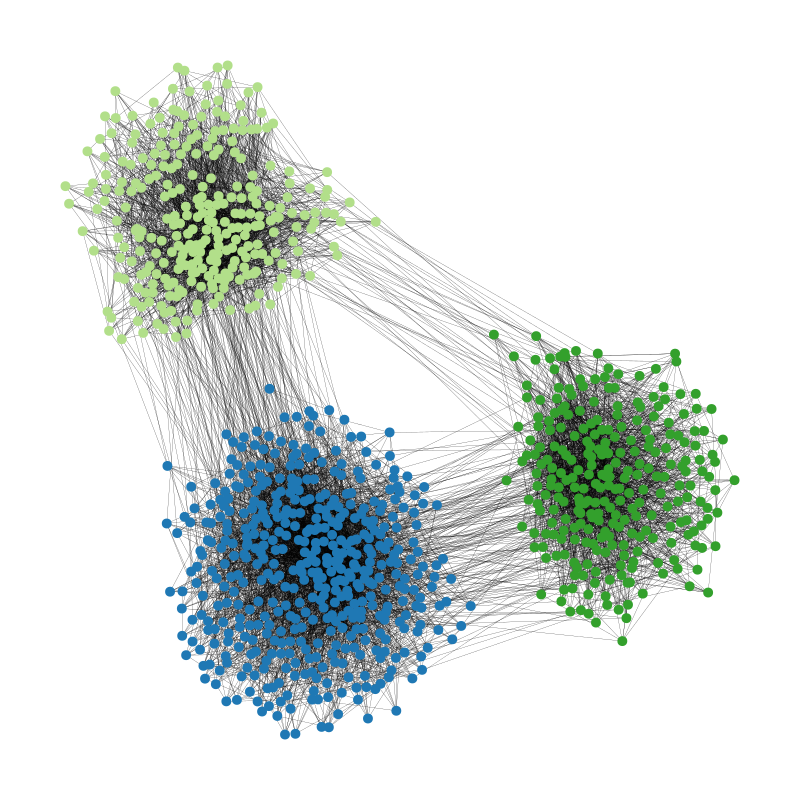}
        \caption{$\mu = \num{5e-2}$}
        \label{fig:lfr graphs 0.05}
    \end{subfigure}
    \begin{subfigure}[b]{0.3\linewidth}
        \includegraphics[width=\linewidth]{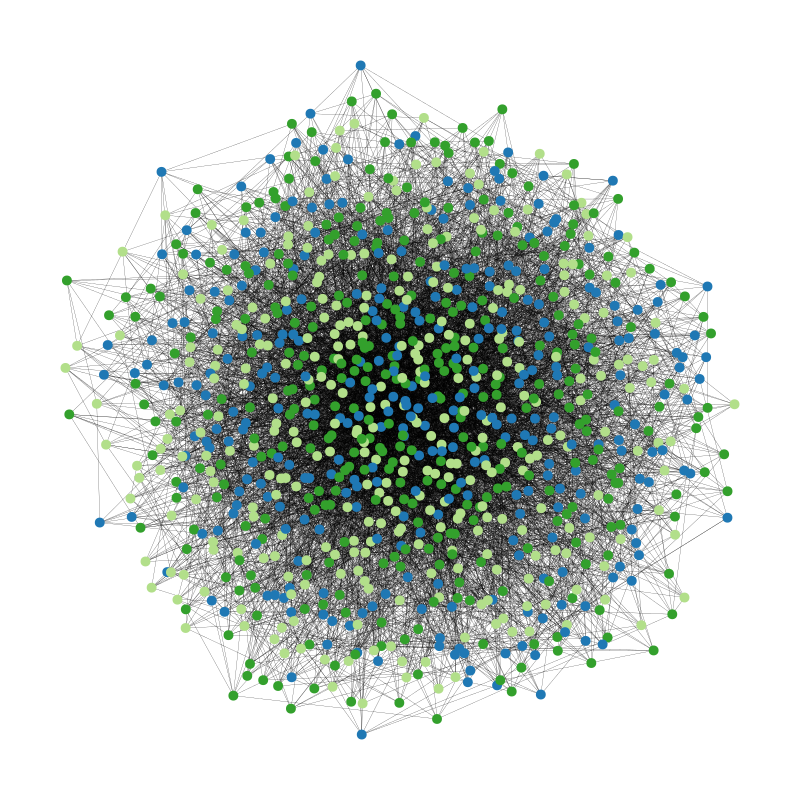}
        \caption{$\mu = \num{5e-1}$}
        \label{fig:lfr graphs 0.525}
    \end{subfigure}
    \caption{Typical LFR graphs for selected values of the mixing parameter $\mu$. We did not include the case $\mu = 1$ as it is visually similar to that of $\mu = 0.5$.}
    \label{fig:lfr graphs}
\end{figure}

We vary the mixing parameter $\mu$ of the LFR algorithm: for $\mu = 0$, individuals are connected to individuals within the same community, whereas for $\mu = 1$, they are connected only to individuals outside their own.
    Some typical network topologies corresponding to $\mu \in \{\num{e-3}, \num{5e-2}, \num{5e-1}, 1\}$ are shown in Figure~\ref{fig:lfr graphs}.

\begin{figure}
  \centering
  \begin{tblr}{
      colspec = {X[c]X[c]X[c]X[c]},
      rowspec = {Q[m]Q[m]Q[m]},
      stretch = 0,
      rowsep = 2pt,
      colsep = 2pt,
    }
    $\mu = \num{e-3}$ &
    $\mu = \num{e-2}$ &
    $\mu = \num{e-1}$ &
    $\mu = \num{5e-1}$ \\
    \includegraphics[width=\linewidth]{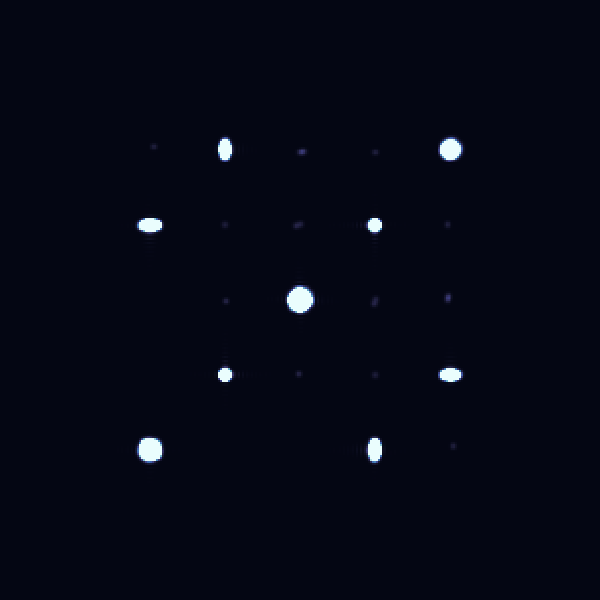} &
    \includegraphics[width=\linewidth]{imgs/lfr/g_mu_0.001.png} &
    \includegraphics[width=\linewidth]{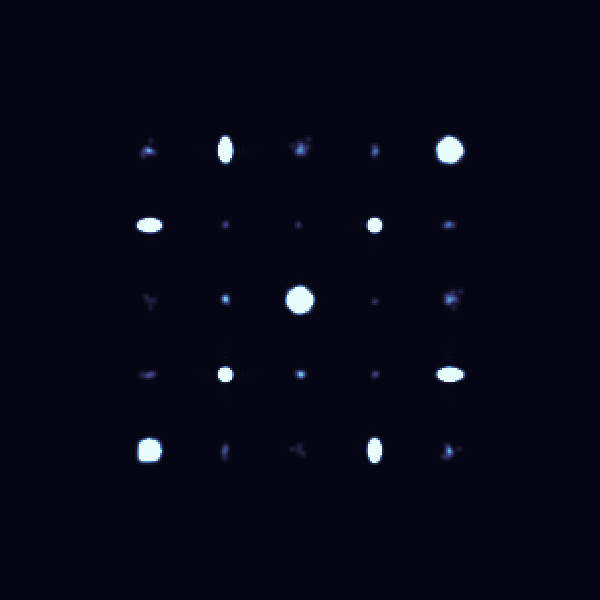} &
    \includegraphics[width=\linewidth]{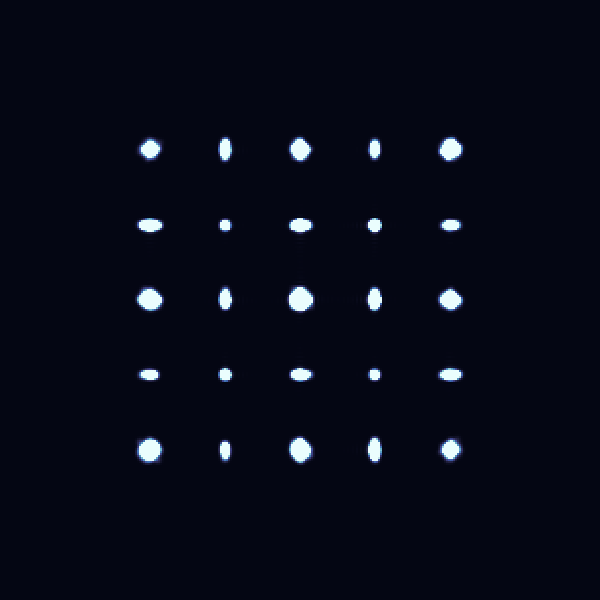} \\
  \end{tblr}
  \caption{Density plots for a typical initial value of $g$. The  lighter the color is, the higher the density is. The value
$\mu$ is the mixing parameter of the LFR graph (see Fig.~\ref{fig:lfr graphs}).
}
  \label{fig:g init}
\end{figure}

    The impact of the value of $\mu$ is reflected in the initial profile for $g$, which is presented in Figure~\ref{fig:g init}. For large $\mu$ (right), all subgroups are connected to each other: for any two $\omega$, $\omega'$ such that $f(\omega)>0$ and $f(\omega')>0$ we have that $g(\omega. \omega')>0$.
    As $\mu$ decreases (going to the left), these links ``fade'': the corresponding $g(\omega, \omega')$ decreases until it reaches $0$ if $(\omega, \omega')$ do not correspond to positive densities in a single subgroup.
  We underline that, in all cases, we consider connected graphs.

\subsection{Results}
\label{sec:numerics results}
Here we present selected results, comparing the discrete and the continuous models on the LFR generated graphs.

The computations were carried out with $N = 1000$ agents for the discrete model and a discretization size of $\Ncont = 303$ for the continuous system. With our choice of parameters, the average degree of nodes is $15$ and the total number of connections $\kappa \sim \num{1.5e4}$,
so that the fill-ratio of the adjacency matrix is roughly 1.5\%. Notice that with these values for $N$ and $\kappa$, the graph is sparse~\cite{diestel_graph_theory_2005}.
We choose $\mathbf{D}(z) := -z$, so that we have the rough estimate $|a[g]| \le 2$, thus the CFL condition associated to the LLF scheme is $\Delta t < \frac{1}{4} \Delta \omega  \simeq 8.10^{-4}$.
This leads us to the choice $\Delta t = 10^{-3}$.

\subsubsection*{Comparison of the convergence to the consensus value $\omega_\infty$}
Because we only consider connected graphs, we expect convergence to consensus in all cases, for both the discrete and continuous models. The consensus opinion is given by $\omega^{micro}_\infty$ and $\omega_\infty^{cont}$, given by \eqref{eq:limit_omega} and \eqref{eq:limit omega continous}. Both coincide, theoretically, as $N\to \infty$.
Because of the sampling required at initialization, there is no reason to expect the two to be equal.

To compare the two models quantitatively, we rather focus on the speed of convergence to consensus
in terms of the $\ell^2$ (resp. $L^2$) distance to $\oinfmicro$ (resp. $\oinfcont$):
$$
  E^\textrm{micro}(t) = \left(\frac{1}{N} \sum_i (\omega_i(t) - \oinfmicro)^2\right)^{\frac{1}{2}}\,, \qquad
  E^\textrm{cont}(t) = \left(\sum_p\int (\omega - \oinfcont)^2 f(t, \omega, p) \xd \omega \right)^{\frac{1}{2}}\,.
$$
In view of Proposition~\ref{prop:expoconv}, we expect exponential convergence of $E^\textrm{micro}$ to zero.

\subsubsection*{Exponential convergence rates and numerical diffusion}
As was already mentioned, the LLF scheme suffers from numerical diffusion, which means that we cannot expect $E^\textrm{cont}$ to vanish but rather to converge to a positive constant, as $t\to \infty$,
something which is confirmed by the plots in Figure~\ref{fig:convergence rates} (left): $E^\textrm{cont}$ decreases but remains roughly above $\num{7e-3}$.

To have a quantitatively meaningful comparison, we estimated the exponential convergence rates for both $E^\textrm{micro}$ and $E^\textrm{cont}$ in the interval $t \in [2, 8]$, in which the effect of the numerical diffusion is only significant for the largest value of $\mu$ (0.5).
There is a transient relaxation phase where the convergence is not exponential (as evidenced by the strictly convex profiles on the left of Figure~\ref{fig:convergence rates}).
This motivates the restriction to $t\in[2,8]$ for the estimation of the convergence rate. In practice, we perform a linear regression between $t$ and $\log E^\textrm{micro}$ (resp. $E^\textrm{cont}$). The fit of this linear regression is very good: for over $500$ sampling points, the average fit error is below $5\%$ for $\mu = 0.5$, and below $1\%$ otherwise.
These convergence rates are shown in Figure~\ref{fig:convergence rates} (right) and the corresponding time interval has been shaded in blue (left).

\begin{figure}
  \begin{tikzpicture}[scale=1]
  \begin{groupplot}[group style={
      group size=2 by 3, horizontal sep=1em,
      y descriptions at=all,
      horizontal sep=0.5cm,
      vertical sep=0.5cm,
    }, grid]
    \nextgroupplot[title={$E^\textrm{micro}, \quad E^\textrm{cont}$}, ymode=log, xmax=8, ymin=0.003, ymax=1, xlabel=$t$]

  \addplot [fill, mark=, color=cyan, opacity=0.5] coordinates { (2.5,0.001) (8,0.001) (8, 2) (2.5,2) } --cycle ;

  \draw[->] (7,0.5) -- (1.2,0.03) node[anchor=north]{$\mu$ incr.};

  \pgfplotstableread{data/sigma_0.012/var_micro.tsv}\micro;
  \pgfplotstableread{data/sigma_0.012/var_mfl_mono.tsv}\contmono;
  \pgfplotstableread{data/sigma_0.012/var_mfl_multi.tsv}\contmulti;

  \foreach \a in {1, 2, 3, 4}{
    \addplot[dotted, thick] table[x index=0, y index=\a] {\micro};
    \addplot[dashed, thick] table[x index=0, y index=\a] {\contmono};
    \addplot[solid, thick] table[x index=0, y index=\a] {\contmulti};
  }
  \nextgroupplot[title={initial convergence rate}, ymode=log, ytick pos=right, legend pos=south east, xmode=log, ymin=0.003, ymax=1.0, xlabel=$\mu$]

    \pgfplotstableread{data/sigma_0.012/rates.tsv}\convrates;

    \addplot[dotted, thick] table[x index=0, y index=1] {\convrates};
    \addlegendentry{$E^\textrm{micro}$}
      \addplot[solid, thick] table[x index=0, y index=2] {\convrates};
    \addlegendentry{$E^\textrm{cont}$ (multi)}
      \addplot[dashed, thick] table[x index=0, y index=3] {\convrates};
    \addlegendentry{$E^\textrm{cont}$ (single)}
\end{groupplot}
\end{tikzpicture}
\caption{Convergence rates as a function of the mixing parameter $\mu$.}
\label{fig:convergence rates}
\end{figure}

\subsubsection*{Observations}

We observe that the rate of convergence increases for increasing values of the connectivity parameter $\mu$ for both the discrete and continuous models. This is not surprising as a larger value of $\mu$ implies increasing connectivity between different groups, particularly, for $\mu=0.5$ one cannot see distinct groups graphically in Figure~\ref{fig:lfr graphs}.

We observe that for small $\mu$ there is a good agreement between the discrete model (DM) and the continuous model with group labeling (CML). In contrast, the continuous model with no group labeling (CMnoL) does not capture the rate of convergence of the microscopic dynamics.
In particular, the estimated convergence rate for CMnoL is practically constant in $\mu$.
This indicates that the refinement of CMnoL to CML is necessary to obtain the right approximation in the limit $N \to \infty$.

For $\mu=0.5$, CML and CMnoL behave similarly, the convergence seems faster than for DM initially, but then slows down until $E^\textrm{cont}$ plateaus. This is also visible on the right of Figure~\ref{fig:convergence rates}, where the estimated convergence rates for CML and CMoL is lower than that of DM. The main cause is the diffusivity of the numerical scheme used for the continuous models, and does not reflect an issue with the models themselves.
Actually, for this value of $\mu$, the error between DM and CMnoL is the smallest. This can be explained, as we mentioned above, by the fact that there are not clear distinct groups in this situation, so adding group labeling does not provide much additional information.

\subsubsection*{Qualitative comparison}

The simulations carried out to compute the convergence graph in Figure~\ref{fig:convergence rates} can be found in Figshare in the form of videos:
\begin{center}
\url{dx.doi.org/10.6084/m9.figshare.28023551}
\end{center}

We can observe qualitatively that the continuous model with group labeling approximates well the discrete dynamics even in the case $\mu=0.5$; and that the continuous model without group labeling clearly fails in all cases except when $\mu=0.5$. We notice also that the case $\mu=0.5$ stands out, as there is no crossing of bumps between the green and blue groups.

\subsubsection*{Interpretation}

Here, we explain why the continuous model with group labeling has a good match with the discrete dynamics, while the model with no group label  fails to do so.

Recall the example of Remark~\ref{rem:validity lemma} and Figure~\ref{fig:snapshots}. Focusing on the histograms of groups $B^\mathrm{in}$ and $Y^\mathrm{in}$ in particular, we see that
they ``cross'' each other and essentially exchange position. There must then be a yellow (resp. blue) individual with index $j_y$ (resp. $j_b$) such that $\omega_{j_y}(t) < \omega_{j_b}(t)$ at $t=0$ and $\omega_{j_y}(t) > \omega_{j_b}(t)$ at $t = 1.2$.

By continuity, there is a time $t^*$ in which $\omega_{j_y}(t^*)=\omega_{j_b}(t^*)$. At this point, it is not possible to distinguish the two individuals at the continuous level. When the opinions overlap, the assumptions of Lemma~\ref{lem:assumption} no longer hold and the network described by $g$ may lose its capacity to represent the original underlying discrete network. However, by labeling the groups, we are able to distinguish between two overlapping opinions and $g$ can now describe with accuracy the underlying network.

To conclude, the numerical simulations indicate that the continuum model with group labeling can be used to investigate the interaction between different communities that initially have non-overlapping opinions with other communities.

\subsection{Implementation and code}
The numerical methods presented above are implemented in Julia, using the following packages:
\texttt{Distributions.jl}~\cite{distributions_jl},
\texttt{Graphs.jl},
\texttt{GraphMakie.jl}\cite{graph_makie},
\texttt{KernelDensity.jl},
\texttt{Makie.jl}\cite{makie},
\texttt{LFRBenchmarkGraphs.jl}\cite{Floros_LFRBenchmarkGraphs_jl},
\texttt{ShiftedArrays.jl},
\texttt{SparseArrays.jl},
\texttt{SpecialFunctions.jl}.
The code is available online\cite{Jankowiak_OpiForm}.

\section{Extensions of the modeling framework to other types of networks}
\label{sec:extensions_model}
The continuum model proposed here is a basic model upon which more complex features can be built. This is one important reason why we proposed this model in contrast to other existing large-particle limit for networks, like graphons. Particularly, there are some important features that can be easily incorporated. We discuss them next.

\subsection{Random creation and destruction of edges (dynamical networks)}
In the spirit of~\cite{degond2016continuum}, one can consider a system where a new edge is created at rate $\lambda$ between  agent $i$ with opinion $\omega$ and agent $j$ with opinion $m$; or that an edge between $\omega$ and $m$ is destroyed at rate $\mu$. These dynamics would result in an additional term on the right-hand-side of equation \eqref{eq:mean-field connectivity} for $g$, of the form:
$$\lambda f(\omega) f(m)-\mu g(\omega,m)\,.$$
This reflects the fact that two opinions $\omega$ and $m$ create a new edge at rate $\lambda$ with a probability given by the proportion of individuals with opinion $\omega$ and $m$, i.e, with a probability $f(\omega)f(m)$; at the same time, an existing edge between opinions $\omega$ and $m$ is destroyed at rate $\mu$ with a probability given by $g(\omega,m)$, i.e. the proportion of edges that exist between opinions $\omega$ and $m$.

One can consider different creation probabilities like, for example, a probability that takes into account the number of intermediate nodes connecting $\omega$ and $m$. In this case we will have a creation term of the form
$$\lambda \left( \int g(\omega, u) g(u,m) du \right)f(\omega) f(m).$$
In another variant, one can consider that edge creation depends on the connectivity of opinion $m$ (where the connectivity corresponds to $\int g(\omega,m)\xd m$), for example, the more connected an opinion is, the more edges are created around that opinion. In this case, we can consider a creation term of the form
$$\lambda \phi[g](\omega,m)\, f(g)f(m),$$
where
$$\phi(g):= H(\max \{\int g(\omega,u) du, \, \int g(m, u) du \})$$
for some increasing function $H$.

One can also consider destruction of edges based on the distance between the opinions. For example, the further opinions are, the easier is to break the edge between them. In this case, one can consider a destruction term of the form
$$-\mu H(|\omega-m|) g(\omega,m),$$
for $H$ increasing.

\subsection{Noise at the particle level}
One can consider a stochastic differential version of our original system of ODEs~\eqref{eq:ODE_eb}:
\begin{align} \label{eq:SDE}
 \xd \omega_i(t)= \frac{1}{\#\mathcal{I}_i} \sum_{j\in \mathcal I_i} \mathbf{D}(\omega_j(t)-\omega_i(t)) \xd t+ \sqrt{2\sigma}dB^i_t,
\end{align}
where $\sigma>0$ is the diffusive constant and $(B^i)_{i=1,\hdots,N}$ are independent Brownian motions. If $\mathbf{D}$ is globally Lipschitz (which is the case of consensus dynamics that we considered here), then we have existence and uniqueness of solutions (pathwise and in law) for \eqref{eq:SDE}.
This would modify (formally) the continuum equations, for the opinion distribution
\begin{align*}
  \partial_t f(\omega) + \partial_\omega \left(a[g](\omega) f(\omega)\right) = \sigma\partial^2_\omega f(\omega)\,, \qquad \omega \in \Omega, \end{align*}
and for the connectivity
$$
  \partial_t g(\omega, m)+ \partial_\omega \left(a[g](\omega) g(\omega,m)\right) + \partial_m \left(a[g](m) g(\omega,m)\right)=\sigma(\partial_\omega^2+\partial_m^2) g(\omega,m)\,, \qquad (\omega, m)\in \Omega^2\,,
$$
where the Laplacian for the equation of $g$ is a consequence of It\^o's Lemma (and adding boundary conditions to have zero flux at the boundary).

Adding diffusion to the dynamics is used to represent self-thinking or changes in opinion due to interactions with the media (which is outside the network), see e.g.~\cite{BS09}.

\subsection{Directional network}
We considered undirected networks, but one can extend this analysis to directional networks. The big difference concerns the definition of $g^N$ in~\eqref{eq:def_gN} which would no longer be symmetric and would be given by
$$g^N(t, \omega, m):= \frac{1}{\kappa}\sum_{i\to j}\delta_{\omega_i(t)}(\omega)\delta_{\omega_j(t)}(m),$$
where $i\to j$ indicates that $i$ influences (is connected to) $j$ (but not the other way around). Then one needs to carry out the steps in Sections~\ref{sec:IBM_nolabels} and~\ref{sec:limit} paying attention to the lack of symmetry of the interactions.

\subsection{Non-weighted interactions}
One can consider the discrete system~\eqref{eq:ODE_eb} without the weights $\# \mathcal I_i$. Also in this setup, one can carry out the analysis as done previously, but the discrete system~\eqref{eq:ODE_eb} needs then to be rescaled. This corresponds to mean-field limits, where the forces (the right-hand-side of~\eqref{eq:ODE_eb}) are rescaled by a multiplying factor $1/N$~\cite{chaintron2021propagation,chaintron2022propagation}. Since we are not working with a mean-field force, but with a network-mediated force, the analog of this would be to rescale the force by a multiplying factor proportional to $1/\kappa$.

\section{Conclusions and perspectives}
\label{sec:conclusions}
In this paper, we propose a framework for the continuum modeling of opinion dynamics on networks, offering an alternative to other approaches, like the ones based on graphons.

We illustrated this framework on a system formed by a fixed undirected network, where individuals undergo pairwise interactions with those individuals with whom they are connected. For this system, we proposed the continuous model~\eqref{eq:kinetic system} as an approximation of the discrete model~\eqref{eq:ODE_eb}. The key for its derivation is to obtain a description of the network that does not use the labels of the individuals. \\
We compared numerically the discrete and continuous models in a setup of consensus dynamics between three communities, and we observed that the continuous model with group labeling gives a good match with the discrete dynamics while the continuous models without group labeling fails in general to do so.
The numerical implementation  required various approximations, and more efficient and accurate methods to simulate the continuum dynamics would be required for more complex setups.

\subsection{Applications}

The goal of the proposed framework is to obtain continuum equations that are amenable to investigating opinion dynamics on social networks.
 For example, as a future perspective, one could investigate reported social effects like epistemic bubbles and echo chambers~\cite{Ng20}. These concepts are used in sociology to explain opinion segregation and polarization via the interaction between communities. We do not dwell here on these concepts, but stress that they have become important in social sciences:
\begin{center}
\begin{minipage}{0.9\linewidth}
``While online echo chambers and epistemic bubbles are commonly found on social media sites like Twitter and Facebook, they are not restricted to such digital platforms. Any online space where information is shared has the potential to become an echo chamber or epistemic bubble. Examples include online forums (e.g., Reddit and Quora), news blogs, video-sharing platforms, online gaming communities, online marketplaces, and review sites.'' (Turner, 2023~\cite{Turner_2023})
\end{minipage}
\end{center}

\subsection{Open problems}

\paragraph{\textbf{The validity of assumption \ref{as:initial_time} as $N\to\infty$.}}

As $N\to\infty$ the conditions for Lemma~\ref{lem:ODE_approx} to hold fall apart since there is no $r>0$ for which~\eqref{eq:def_r} holds as $N\to \infty$. The hope, though, is that the continuum equations provide a good approximation of the discrete ones. For the moment, we can test this only numerically. It will be the issue of future investigation to find a correct analytic measurement to quantify how good is the continuum approximation.

\paragraph{\textbf{Conditions for the limiting network to exist and rigorous mean-field limit.}}

A fundamental question is how a network should scale as \( N \to \infty \) for the continuum equations to make sense. In other words, we require that the initial edge distribution \( g^N|_{t=0} \) admits a well-defined limit as \( N \to \infty \). What assumptions on \( g^N|_{t=0} \) are necessary to guarantee the existence of this limit?

Even if this limit exists, rigorously proving the convergence of \( f^N \) and \( g^N \) to their respective continuum counterparts as \( N \to \infty \) (Assumption~\ref{as:limit}) remains a significant challenge~\cite{chaintron2021propagation,chaintron2022propagation}. Ideally, one would also like to derive an analytical estimate of the convergence rate of the discrete system to the continuum system for large values of \( N \), similar to what is achieved in classical mean-field limits~\cite{carmona2016lectures,sznitman1991topics}.

\paragraph{\textbf{Numerical methods.}} As mentioned in Section~\ref{sec:refinements}, refinements of the basic framework may be needed to capture well the network structure. These refinements require adding more information in the system, which increases its computational cost. A future problem will be how to find efficient numerical methods for some of these refinements.

\section*{Acknowledgments}
The authors thank Pierre Degond (CNRS, University of Toulouse) and Antoine Diez (University of Kyoto) for useful discussions and tips. GF has been partially financed by the Austrian Science Fund (FWF) project \href{https://www.fwf.ac.at/en/research-radar/10.55776/F65}{10.55776/F65}, by INdAM | Gruppo Nazionale per l'Analisi Matematica, la Probabilità e le loro Applicazioni (INdAM/GNAMPA) no.~E55F22000270001 and no.~E53C22001930001. GF also wants to thank the Department of Applied Mathematics and Mathematical Physics of the Urgench State University. The work of SMA was funded in part  by the Vienna Science  and  Technology  Fund  (WWTF)  \href{https://www.wwtf.at/funding/programmes/vrg/VRG17-014/index.php?lang=EN}{10.47379/VRG17014} and in part by the Austrian Science Fund (FWF) project \href{https://www.fwf.ac.at/en/research-radar/10.55776/F65}{10.55776/F65}.

\section*{Conflict of Interest Statement}
The authors have no conflicts to disclose.

\begin{appendix}
  \section{Exponential convergence to consensus in the linear case.}
  \label{sec:app:expoconv}
  Here we sketch a proof of Proposition~\ref{prop:expoconv}.
  In the case $\mathbf{D}(z) = -z$, our system of equations~\eqref{eq:ODE_eb} is a linear ODE system of the form
\[
\begin{cases}
  \frac{\xd u(t)}{\xd t} = Qu(t) & t>0\,,
  \\
  u(0) = u_0 & t=0\,,
\end{cases}
\]
where $u=(\omega_1,\hdots, \omega_N)$ and the constant $N\times N$ matrix $Q$ is given by
\[ Q_{ij} = \begin{cases}
  -1 & \text{ if }\quad  i = j, \\
  (\# \mathcal I_i)^{-1} & \text{ if }\quad i \sim j. \\
  0 & \text{ otherwise}\,.
\end{cases}
\]
One has $Q = M^{-1}A - I$, where $I$ is the identity matrix and where
$M := \diag(\#\mathcal{I}_i)$ is the diagonal matrix whose $i^\text{th}$ diagonal element is $\# \mathcal{I}_i$.
With this definition, the matrix $L := -M Q = -Q^T M$ is called the \emph{Laplacian matrix} of the graph in spectral graph theory:
\[ L_{ij} = \begin{cases}
  \#\mathcal{I}_i & \text{ if }\quad  i = j, \\
  -1 & \text{ if }\quad i \sim j, \\
  0 & \text{ otherwise}\,.
\end{cases}
\]
The matrix $L$, unlike $Q$, is symmetric.
It holds that $L_{ii} = \sum_{j\neq i}|L_{ij}|$ and $Q_{ii} = -\sum_{j \neq i}|Q_{ij}|$.
As a consequence, through Gershgorin's Disc Theorem~\cite{gerschgorin_ev_1931}, we know that their respective spectra, denoted by $\sigma(L)$ and $\sigma(Q)$, satisfy
\[
  \sigma(L) \in [0, 2 \max \# \mathcal{I}_i], \quad \real(\sigma(Q)) \in [-2, 0]\,,
\]
where $\real$ denotes the real part. For $L$, the multiplicity of the eigenvalue $0$ is given by the number of connected components of the graph. In our case the graph is connected,
so $0$ is a simple eigenvalue. Since $\ker L = \ker Q$, $0$ is also a simple eigenvalue of $Q$.
The corresponding eigenspace ---the kernel--- is given by $\Span(\mathbf{1}_n)$, where $\mathbf{1}_n = (1,\dots,1)$, i.e. vectors whose entries are equal.

We can write $u$ using a decomposition on the kernel and its orthogonal:
\[u = \tilde{u} + \oinfmicro \mathbf{1}_n := (u - \oinfmicro \mathbf{1}_n) + \oinfmicro\mathbf{1}_n \,,\]
by construction, $\tilde u(t) \in \ker(Q)^\perp$ and $\oinfmicro \mathbf{1}_n \in \ker(Q)$ for all $t$.
It then holds
\[ \frac{d \tilde{u}}{dt} = \frac{d u}{dt} \,,\]
again by construction.
We can now show the exponential convergence of a weighted $\ell^2$-norm of $\tilde{u}$. We introduce the inner product $\langle \cdot, \cdot \rangle_M := \langle \cdot, M \cdot \rangle$,
along with the associated norm, $\| \cdot \|_M$.
Because $L = -MQ$ is symmetric, $Q$ is self-adjoint with this inner product. We have
$$
  \frac{1}{2} \frac{d}{dt} \| \tilde u\|_M^2
  = \langle \tilde{u}, Q \tilde{u}\rangle_M = -\langle\tilde{u}, L \tilde{u}\rangle
  \le -\min_{\lambda \in \sigma(L) \setminus \{0\}} \lambda\ \; \|\tilde{u}\|^2
   \le -\left(\min_{\lambda \in \sigma(L) \setminus \{0\}}\lambda \; \min_i \# \mathcal{I}_i\right) \|\tilde{u}\|_M^2 \,.
$$
The first inequality follows from the fact that $\tilde{u} \in \ker(L)^\perp$. The second from the fact that $\# \mathcal{I}_i \ge 1$ by assumption. Gronwall's Lemma yields \emph{exponential convergence}
of $\|\tilde{u}\|_M$:
\begin{equation}
  \|\tilde{u}(t)\|_M \le \exp\left(- \alpha t\right) \|\tilde{u}_0\|_M\,,
\end{equation}
where $\alpha = \min_{\lambda \in \sigma(L) \setminus \{0\}} \min_{i} \# \mathcal{I}_i$ is the exponential convergence \emph{rate}.
By the equivalence of norms, this also holds for the Euclidean norm $\|\cdot\|$ up to a multiplicative constant $C$.

  \section{Step-size for the discrete model}
  \label{sec:app:stepsize}
  In this Section, we investigate which condition on the step-size $\Delta t$ is sufficient to get rid of the clamping function in \eqref{eq:micro:discretized}.
  More precisely, we prove

  \begin{proposition}
    Consider a connected graph of size $N$, with opinions $(\omega_i)_{1 \le i \le N} \subset \Omega$ and $(\mathcal{I}_i)_{1 \le i \le N}$ defined as in Section~\ref{sec:iba}.
    Let $\mathbf{D} : \mathbb{R} \to \mathbb{R}$ satisfying Assumption~\ref{ass:debate operator}. For $1 \le i \le N$, define $\psi_i$
    \[ \psi_i := \omega_i + \Delta t \frac{1}{\# \mathcal{I}_i} \sum_{j\in \mathcal{I}_i} \mathbf{D}(\omega_i - \omega_j)\,.\]
    If $\Delta t \le \| \partial_z\mathbf{D}\|_{L^\infty}^{-1} = \| \partial_z^2 W\|_{L^\infty}^{-1}$, then $(\psi_i)_{1 \le i \le N} \subset \conv(\{\omega_i\} \cup \{\omega_j : j \in \mathcal{I}_i \})$,
    where $\conv(A)$ denotes the convex hull of the set $A$.

    In particular, we have $(\psi_i)_{1 \le i \le N} \subset \Omega$.
  \end{proposition}

  \begin{proof}
    We define $[c_i^{\min}, c_i^{\max}] := \conv(\{\omega_i\} \cup \{\omega_j : j \in \mathcal{I}_i \})$,
    so that we must show:
    \begin{equation*}
      c_i^{\min} - \omega_i \le \psi_i - \omega_i \le c_i^{\max} - \omega_i\,, \quad \forall 1 \le i \le N\,.
    \end{equation*}
    We only deal with the lower bound (the left inequality), the upper bound can be shown in an identical way.
    First define $\underline{\mathbf{D}} : z \mapsto - \|\partial_z^2 W\|_{L^\infty} \max(0, z)$, which is non-increasing. We have
    $$
      \mathbf{D}(\omega_i - \omega_j) = \int_0^{\omega_i - \omega_j} \partial_z \mathbf{D}(z) \; dz                                      = -\int_0^{\omega_i - \omega_j} \partial_z^2 \mathbf{W}(z) \; dz                                       \ge \underline{\mathbf{D}}(\omega_i - \omega_j) \ge \underline{\mathbf{D}}(\omega_i - c_i^{\min})\,,
    $$
    using the fact that $\partial_z^2 W > 0$. Summing over $j \in \mathcal{I}_i$ and dividing by $\#\mathcal{I}_i$ we get
    \[ \psi_i - \omega_i \ge \Delta t \, \underline{\mathbf{D}}(\omega_i - c_i^{\min}) \,. \]

    On the one hand, if $c_i^{\min} < \omega_i$, then $\underline{\mathbf{D}}(\omega_i - c_i^{\min}) = \| \partial_z^2 W\|_{L^\infty} (c_i^{\min} - \omega_i) \le 0$ and it follows
    \[ c_i^{\min} - \omega_i \le \Delta t \, \|\partial_z^2 W\|_{L^\infty} (c_i^{\min} - \omega_i) \le \psi_i - \omega_i \,.\]

    On the other hand, if $c_i^{\min} - \omega_i = 0$, i.e. $\omega_i - \omega_j \le 0$ for $1\le j \le N$, then $0 \le \mathbf{D}(\omega_i - \omega_j)$ and $c_i^{\min} - \omega_i = 0 \le \psi_i - \omega_i$, which completes the proof.
  \end{proof}

\end{appendix}

\bibliographystyle{abbrv}
\bibliography{ref}

\end{document}